\documentclass[a4paper,11pt,twoside]{article}
\usepackage{pst-fill,pst-grad}
\usepackage{textcomp}
\usepackage[english]{babel}
\usepackage[utf8x]{inputenc}
\usepackage{graphicx}
\usepackage{multicol}
\usepackage{float}
\usepackage{fancyhdr}
\usepackage[matrix,arrow,curve]{xy}
\usepackage{pstricks} 
\usepackage{amsmath,amsfonts,verbatim,afterpage,euscript,mathrsfs,amssymb}
\usepackage{amsthm}
\usepackage{array}
\usepackage{dsfont}
\usepackage{hyperref}

\newcommand{\mysection}{\setcounter{equation}{0} \section}

\def \vu{\vec{u}}
\def \vw{\vec{\omega}}
\def \va{\vec{a}}
\def \vb{\vec{b}}
\def \vc{\vec{c}}
\def \vf{\vec{f}}
\def \vg{\vec{g}}

\def \vA{\vec{A}}
\def \vB{\vec{B}}
\def \vR{\vec{R}}
\def \vU{\vec{U}}
\def \vMU{\vec{\mathcal{U}}}

\def \vW{\vec{W}}

\def \vn{\vec{\nabla}}
\def \vphi{\vec{\varphi}}
\def \rot{\vec{\nabla}\wedge}

\def \R{\mathbb{R}^{3}}
\def \M{\mathcal{M}_{t,x}}

\newtheorem{propo}{Proposition}[section]
\newtheorem{lem}{Lemma}[section]
\newtheorem{theo}{Theorem}[section]
\newtheorem{coro}{Corollary}[section]
\newtheorem{rema}{Remark}[section]

\title{\bf A crypto-regularity result for the\\ micropolar fluids equations}
\usepackage{authblk}
\author{Diego Chamorro\footnote{\emph{diego.chamorro@univ-evry.fr}} }
\author{David Llerena\footnote{\emph{david.llerena@univ-evry.fr}} }
\affil{\footnotesize LaMME, Univ. Evry, CNRS, Universit\'e Paris-Saclay, 91025, Evry, France.}

\begin{document}
\sloppy
\maketitle
\begin{scriptsize}
\abstract{In the analysis of PDEs, regularity of often measured in terms of Sobolev, H\"older, Besov or Lipschitz spaces, etc. However, sometimes a gain of regularity can also be expressed just in terms of Lebesgue spaces, by passing from a singular setting to a less singular one. In this article we will obtain a gain of integrability for weak solutions of the micropolar fluid equations using as general framework Morrey spaces, which is a very useful language to study regularity in PDEs. An interesting point is that the two variables of the micropolar fluid equations can be studied separately.}\\

\end{scriptsize}

%\tableofcontents
%%%%%%%%%%%%%%%%%%%%%%%%%%%%%%%%%%%%%%%%%%%%%%%%%%%
\mysection{Introduction}
In this article we study some integrability results for weak solutions of the perturbed 3D micropolar fluid equations given by the following system:
\begin{equation}\label{MicropolarEquations}
\begin{cases}
\partial_t \vu  = \Delta \vu  -(\vu \cdot \vn)\vu-\vn p +\frac{1}{2}\rot\vw+ \vf + (\va\cdot \vn)\vu+(\vu\cdot\vn)\va, \qquad div(\vu)=div(\va)=0,\\[3mm]
\partial_t \vw = \Delta \vw +\vn div(\vw)-\vw -(\vu \cdot \vn)\vw+\frac{1}{2}\rot\vu, \qquad\qquad\qquad\quad div(\vf)=0,\\[3mm]
\vu (0,x)=\vu_0(x),\;\; div(\vu _0)=0, \quad \vw(0,x)= \vw_0(x),\; x\in \R.
\end{cases}
 \end{equation}
In the previous equations $\vu$ is the fluid velocity field, $\vw$ is the field of microrotation representing the angular velocity of the rotation of the fluid particles, $p$ is the scalar pressure, the quantity $\vf$ represents some given divergence free external force and $\vu_0$ and $\vw_0$ are the initial data. The perturbation $\va\in L^6_tL^6_x$ which appears in the first equation above in the term $ (\va\cdot \vn)\vu+(\vu\cdot\vn)\va$ is a given divergence free vector field and the presence of this particular type of perturbation is inspired by some applications (see \cite{BP2} and \cite[Section 12.6]{PGLR1}) and it is linked to a previous work \cite{ChLLKukaMMP}. Of course, if $\vw=\va=0$ this system reduces to the classical Navier-Stokes equations.\\

Micropolar equations were first introduced in 1966 by Eringen \cite{Eringen} and they have been recently studied by many authors, see \emph{e.g.} \cite{ChLLKukaMMP}, \cite{Galdi}, \cite{GuWang}, \cite{LoMelo} and the references there in. One very interesting feature of the micropolar fluid system is that the variable $\vw$ is not a divergence free vector field and this fact makes its study quite different from other systems of PDEs based on the Navier-Stokes equations (such as the Magneto-Hydrodyamic equations for example). \\

The main result of this article (see Theorem \ref{Theo1} below) states that, if we have some local control of the variable $\vu$ in terms of parabolic Morrey spaces (see formula (\ref{DefMorreyparabolico} below for a definition of these spaces), then we can deduce a local control for the variables $\vu$ and $\vw$ in terms of Lebesgue spaces. 
Let us remark that Morrey spaces are generalization of Lebesgue spaces and they can contain quite singular objects (see the books \cite{Adams, PGLR1, Triebel} for a detailed treatment of Morrey spaces) thus, this gain of integrability can be seen as a gain of \emph{crypto}-regularity as we pass from a singular framework to a less singular framework (and this motivates the title of the article). To be more precise, recall that in the setting of the 3D Navier-Stokes equations it is known that if we have some suitable local integrability information then, following the work of Serrin \cite{Serrin1} in the case $\frac{2}{p}+\frac{3}{q}<1$ with $q>3$ and Struwe \cite{Struwe} and Takahashi \cite{Takahashi} in the case $\frac{2}{p}+\frac{3}{q}=1$ with $q>3$, it is possible to deduce a (usual) gain of regularity in terms of Sobolev spaces. These results were generalized by O'Leary \cite{OLeary}  to the framework of Morrey spaces and it is shown there that a local gain of integrability implies, by the Serrin local regularity theory (see \cite[Chapter 13]{PGLR1}), a gain of regularity. Note that, to the best of our knowledge, in the setting of the perturbed equation (\ref{MicropolarEquations}) all the different results that constitute the Serrin theory are not available, thus we will follow here the general spirit of this theory, but with some modifications. Indeed, due to the structure of the system (\ref{MicropolarEquations}) we will perform our study in two separated steps: first we will study the variable $\vu$ in order to obtain a gain of integrability on $\vu$ and only then we will use this result to study the variable $\vw$ and doing so, we will obtain a gain of integrability for the two variables, note that at each step we will also need to deal with the perturbation term and the external forces. Finally, we will show how to obtain a classical H\"older regularity result for the couple $(\vu, \vw)$. However, we will not look for a Serrin type regularity theory and instead we will use the more sophisticated theory of Caffarelli-Kohn-Nirenberg which was displayed in our previous article \cite{ChLLKukaMMP}: doing so we will mix these two theories in order to obtain a different and new regularity framework for the system (\ref{MicropolarEquations}).\\

Let us recall now the definition of parabolic Morrey spaces. For $1< p\leq q<+\infty$, parabolic Morrey spaces $\mathcal{M}_{t,x}^{p,q}$ are defined as the set of measurable functions $\vphi:\mathbb{R}\times\R\longrightarrow \R$ that belong to the space $(L^p_{t,x})_{loc}$ such that $\|\vphi\|_{\mathcal{M}_{t,x}^{p,q}}<+\infty$ where
\begin{equation}\label{DefMorreyparabolico}
\|\vphi\|_{\mathcal{M}_{t,x}^{p,q}}=\underset{x_{0}\in \R, t_{0}\in \mathbb{R}, r>0}{\sup}\left(\frac{1}{r^{5(1-\frac{p}{q})}}\int_{|t-t_{0}|<r^{2}}\int_{B(x_{0},r)}|\vphi(t,x)|^{p}dxdt\right)^{\frac{1}{p}}.
\end{equation}
Note that we have $\mathcal{M}_{t,x}^{p,p}=L_{t,x}^p$. See \cite{Adams} for more details on these spaces or \cite{PGLR1} for a general theory concerning the Morrey spaces and Hölder continuity applied to the analysis of PDEs from fluid mechanics. See also Section \ref{Secc_Morrey} below for a list of useful facts concerning these spaces.\\

As we are interested in the study of  \emph{local} properties of a weak solution $(\vu,p,\vw)$ to the micropolar fluid equations \eqref{MicropolarEquations}, we will fix once and for all an open parabolic ball  $Q \subset]0,+\infty[\times \R$ of the form 
\begin{equation}\label{SetQ}
Q = ]a,b[\times B(x_0, r), \quad \mbox{with}\quad 0<a<b<+\infty,\;\; x_0\in \R \;\; \mbox{and} \;\; 0<r<+\infty,
\end{equation}
note in particular that ($0,0)\notin Q$. We are now ready to state our main result in which we obtain a gain of integrability for $\vu$ and for $\vw$.
%%%%%%%%%%%%%%%%%%%%%%%%%%%%%%%%%%%%%%%%%%%%%%%%%%%
\begin{theo}\label{Theo1}
Let $(\vu,p, \vw)$ be a weak solution over $Q$ of the perturbed micropolar fluid equations (\ref{MicropolarEquations}) where $Q$ is the parabolic ball given in (\ref{SetQ}). Assume that $\vu,\vw\in L_t^\infty L_x^2 \cap L_t^2 \dot H_x ^1(Q)$, $p\in \mathcal{D}'_{t,x}(Q)$, $\vf\in L^2_t\dot{H}^1_x(Q)$ and $\va\in L^6_tL^6_x(Q)$.
If moreover we have the following local hypothesis
\begin{equation}\label{Hypothese1}
\mathds{1}_Q\vu \in \M^{p_0,q_0}(\mathbb{R}\times\R)\quad with\;\;
2<p_0\le q_0,\; 5<q_0\le 6,
\end{equation}
then
\begin{itemize}
\item[1)] for a parabolic ball $Q_1=]a_1,b_1[\times B(x_0, r_1)$, with  $a<a_1< b_1<b$ and $0<r_1<r$ (we thus obtain the inclusion $Q_1\subset Q$), we have that 
\begin{equation}\label{Conclusion1}
\mathds{1}_{Q_1}\vu \in L_{t,x}^{q_0}(\mathbb{R}\times\R),\qquad 5<q_0\le 6.
\end{equation}
\item[2)] for a parabolic ball $Q_2=]a_2,b_2[\times B(x_0, r_2)$, with  $a_1<a_2< b_2<b_1$ and $0<r_2<r_1$ (note the inclusion of parabolic balls $Q_2\subset Q_1$), we have that 
\begin{equation}\label{Conclusion2}
\mathds{1}_{Q_2}\vw \in L_{t,x}^{q_0}(\mathbb{R}\times\R),
\end{equation}
for some $5<q_0\leq 6$.
\end{itemize}
\end{theo}
%%%%%%%%%%%%%%%%%%%%%%%%%%%%%%%%%%%%%%%%%%%%%%%%%%%
\noindent Some remarks are in order here. Note first that we only impose the Morrey-type condition (\ref{Hypothese1}) to the variable $\vu$ and we do not assume any further assumption over $\vw$. Indeed, with the sole information $\vw\in L_t^\infty L_x^2 \cap L_t^2 \dot H_x ^1(Q)$ we can perform our computations for the variable $\vu$ and we can obtain for this variable the local gain of integrability given in (\ref{Conclusion1}). Now, once we have this information over $\vu$ at hand, we can focus on the variable $\vw$ and we can prove (\ref{Conclusion2}) and this two-step procedure shows that, when studying the gain of integrability problem, the velocity field $\vu$ ``dominates'' the angular velocity $\vw$. 

Remark now that the pressure $p$ is a very general object as we only ask $p\in \mathcal{D}'$, this is a general feature of the Serrin local regularity theory and the usual trick to get rid of the pressure consists in applying the Curl operator to the first equation in (\ref{MicropolarEquations}). However, if we want to mix our results with other regularity theories (such as the Caffarelli-Kohn-Nirenberg criterion), then some information over the pressure will be needed, see Corollary \ref{Coro_Application} below. To end the remarks, it is worth noting here that the upper bound for the parameter $q_0$ in (\ref{Hypothese1}), \emph{i.e.} the range $5<q_0\leq 6$, is mainly technical and it is related to the information available on the perturbation term $\va\in L^6_tL^6_x$. Let us mention that this control over $\va$ is very useful to obtain existence of Leray-type solution to the problem (\ref{MicropolarEquations}) and we do not claim any optimality over these parameters.\\

Let us now explain how the previous crypto-regularity result (Theorem \ref{Theo1}) can help to deduce a \emph{real} gain of regularity in terms of usual spaces. The following result heavily relies in our previous work done in \cite{ChLLKukaMMP}.
%%%%%%%%%%%%%%%%%%%%%%%%%%%%%%%%%%%%%%%%%%%%%%%%%%%
\begin{coro}\label{Coro_Application}
 Let $(\vu,p, \vw)$ be a weak solution over $Q$ of the perturbed micropolar fluid equations (\ref{MicropolarEquations}) where $Q$ is the parabolic ball given in (\ref{SetQ}). Assume that $\vu,\vw\in L_t^\infty L_x^2 \cap L_t^2 \dot H_x ^1(Q)$, $p\in \mathcal{D}'_{t,x}(Q)$, $\vf\in L^2_t\dot{H}^1_x(Q)$ and $\va\in L^6_tL^6_x(Q)$. Assume moreover the following points:
\begin{itemize}
\item[1)] $\mathds{1}_Q\vu \in \M^{p_0,q_0}(\mathbb{R}\times\R)$ with $2<p_0\le q_0,\; 5<q_0\le 6$,
\item[2)] The pressure $p$ belongs to the space $L^{\frac{3}{2}}_{t,x}(Q)\cap L^{\frac{5}{2}}_{t}L^1_x(Q)$,
\item[3)] The external force satisfies $\mathds{1}_Q\vf\in \M^{\frac{10}{7}, \tau}(\mathbb{R}\times \R)$ for some $\tau>\frac{5}{2-\alpha}$ with $0<\alpha\ll1$,
\item[4)] The weak solution $(\vu, p, \vw)$ is suitable in the sense of the Caffarelli-Kohn-Nirenberg partial regularity theory (see the Definition 1.1 of \cite{ChLLKukaMMP}, see also the Section 13.8 of the book \cite{PGLR1}).
\item[5)] There exists a positive constant $\epsilon^*$ which depends only on $\tau$ such that, if for some $(t_0,x_0)\in Q$, we have
\begin{equation*}
\limsup_{r\to 0}\frac{1}{r}\iint _{]t_0-r^2, t_0+r^2[\times B(x_0,r)}|\vn \otimes \vu|^2+|\vn \otimes \vw|^2dx ds<\epsilon^*,
\end{equation*}
then $(\vu,\vw)$ is H\"older regular (in the time and space variables) of exponent $\alpha$ in a neighborhood of $(t_0,x_0)$ for some small $\alpha$ in the interval $0<\alpha\ll 1$.
\end{itemize}
\end{coro}
Note in particular that, to obtain a H\"older regularity result for the variables $\vu$ and $\vw$, we only need to impose an integrability hypothesis on $\vu$.\\

\noindent {\bf Proof of the Corollary \ref{Coro_Application}.} Remark that the first hypothesis here corresponds to the condition (\ref{Hypothese1}) above and thus, applying Theorem \ref{Theo1} we obtain the conclusion (\ref{Conclusion2}), \emph{i.e.} $\mathds{1}_{Q_2}\vw \in L_{t,x}^{q_0}(\mathbb{R}\times\R)$ for some $5<q_0\leq 6$. Now, if we fix $0<\alpha\ll 1$ small enough in order to verify the condition $\frac{5}{q_0}<1-2\alpha$, we actually satisfy a local condition of integrability for $\vw$ needed in order to apply the Caffarelli-Kohn-Nirenberg regularity theory, see Remark 4.2 of \cite{ChLLKukaMMP}. Thus, following this theory displayed in the article \cite{ChLLKukaMMP} we obtain the wished H\"older regularity result for the variables $\vu$ and $\vw$. \hfill $\blacksquare$\\

The plan of the article is the following: in Section \ref{Secc_Morrey} we recall some useful facts related to the Morrey spaces. In Section \ref{Secc_GainU} we study the gain of integrability for the variable $\vu$ while in Section \ref{Secc_GainOmega} we obtain the gain of integrability for the variable $\vw$.
%%%%%%%%%%%%%%%%%%%%%%%%%%%%%%%%%%%%%%%%%%%%%%%%%%%
\mysection{Useful results related to parabolic Morrey spaces}\label{Secc_Morrey}
In this section we state without proof some results on parabolic Morrey spaces and we refer to the books \cite{Adams, PGLR1} for a proof of these facts and for a more detailed study of these functional spaces. We start recalling the parabolic framework given by the homogeneous space $(\mathbb{R}\times \R, d, \mu)$ where $d$ is the parabolic quasi-distance given by 
\begin{equation}\label{Def_QuasiDistance}
d\big((t,x), (s,y)\big)=|t-s|^{\frac{1}{2}}+|x-y|,
\end{equation}
and where $\mu$ is the usual Lebesgue measure $d\mu=dtdx$. Remark that the homogeneous dimension is $N=5$. See \cite{Folland} for more details concerning the general theory of homogeneous spaces. 
%%%%%%%%%%%%%%%%%%%%%%%%%%%%%%%%%%%%%%%%%%%%%%%%%%%
\begin{lem}[H\"older inequalities]\label{lemma_Product}
\begin{itemize}
\item[]
\item[1)]If $\vf, \vg:\mathbb{R} \times \R\longrightarrow \R$ are two functions such that $\vf\in \M^{p,q} (\mathbb{R} \times \R)$ and $\vg\in L^{\infty}_{t,x} (\mathbb{R} \times \R)$, then for all $1\leq p\leq q<+\infty$ we have
$\|\vf\cdot\vg\|_{\M^{p,q}}\leq  C\|\vf\|_{\M^{p, q}} \|\vg\|_{L^{\infty}_{t,x}}$.
\item[2)]If $\vf, \vg:\mathbb{R} \times \R\longrightarrow \R$ are two functions that belong to the space $\M^{p,q} (\mathbb{R} \times \R)$ then we have the inequality $\|\vf\cdot\vg\|_{\M^{\frac{p}{2}, \frac{q}{2}}}\leq  C\|\vf\|_{\M^{p, q}} \|\vg\|_{\M^{p, q}}$.
\item[3)] More generally, let $1\leq p_0 \leq q_0 <+\infty$, $1\leq p_1\leq q_1<+\infty$ and $1\leq p_2\leq q_2<+\infty$. If $\tfrac{1}{p_1}+\tfrac{1}{p_2}\leq \frac{1}{p_0}$ and $\tfrac{1}{q_1}+\tfrac{1}{q_2}=\tfrac{1}{q_0}$, then for two functions $\vf, \vg:\mathbb{R} \times \R\longrightarrow \R$ such that $\vf\in \mathcal{M}^{p_1, q_1}_{t,x} (\mathbb{R} \times \R)$ and $\vg\in \mathcal{M}^{p_2, q_2}_{t,x} (\mathbb{R} \times \R)$, we have the following H\"older inequality in Morrey spaces 
$$\|\vf\cdot \vg\|_{\mathcal{M}^{p_0, q_0}_{t,x}}\leq \|\vf\|_{\mathcal{M}^{p_1, q_1}_{t,x}}\|\vg\|_{\mathcal{M}^{p_2, q_2}_{t,x}}.$$
\end{itemize}
\end{lem} 
%%%%%%%%%%%%%%%%%%%%%%%%%%%%%%%%%%%%%%%%%%%%%%%%%%%
\begin{lem}[Localization]\label{lemma_locindi}
Let $Q$ be a bounded set of $\mathbb{R} \times \R$ of the form given in  (\ref{SetQ}). If we have $1\leq p_0 \leq p_1$ and $1\leq p_0\leq q_0 \leq q_1<+\infty$ and if the function $\vf:\mathbb{R} \times \R\longrightarrow \R$  belongs to the space $\M^{p_1,q_1} (\mathbb{R} \times \R)$ then we have the following localization property 
$$\|\mathds{1}_{Q}\vf\|_{\M^{p_0, q_0}} \leq C\|\mathds{1}_{Q}\vf\|_{\M^{p_1,q_1}}\leq C\|\vf\|_{\M^{p_1,q_1}}.$$
\end{lem}  
%%%%%%%%%%%%%%%%%%%%%%%%%%%%%%%%%%%%%%%%%%%%%%%%%%%
\noindent Now, for $0<\mathfrak{a}<5$ we can define the (parabolic) Riesz potential $\mathcal{I}_{\mathfrak{a}}$ of a locally integrable function
$\vec f: \mathbb{R}\times\mathbb{R}^3\longrightarrow \mathbb{R}^3$ by
\begin{equation}\label{DefinitionPotentielRiesz6}
\mathcal{I}_{\mathfrak{a}}(\vf)(t,x)=\int_{\mathbb{R}}\int_{\mathbb{R}^3}
\frac{1}{(|t-s|^{\frac{1}{2}}+|x-y|)^{5-\mathfrak{a}}}\vec{f}(s,y)dyds.
\end{equation}
Then, we have the following property
%%%%%%%%%%%%%%%%%%%%%%%%%%%%%%%%%%%%%%%%%%%%%%%%%%%
\begin{lem}[Adams-Hedberg inequality]\label{Lemme_Hed}
If $0<\mathfrak{a}<\frac{5}{q}$, $1<p\le q<+\infty$ and $\vf\in \M^{p,q} (\mathbb{R} \times \R)$,
then for $\nu=1-\frac{\mathfrak{a}q}{5}$ we have the following boundedness property in Morrey spaces:
\begin{equation*}
\|\mathcal{I}_{\mathfrak{a}}(\vf)\|_{\M^{\frac{p}{\nu},\frac{q}{\nu}}}\le\|\vf\|_{\M^{p,q}}. 
\end{equation*}
\end{lem}
\noindent We finish this section with two simple corollaries of the previous lemma obtained in \cite{LocalMHD} and \cite{PGLR1} which will be helpful in our computations
%%%%%%%%%%%%%%%%%%%%%%%%%%%%%%%%%%%%%%%%%%%%%%%%%%%
\begin{coro}[Estimates for $\mathcal{I}_{1}$]\label{Morrey_Coro1}
Let $Q$ be a bounded set of $\mathbb{R} \times \R$ of the form given in  (\ref{SetQ}). If we have $2< p \le q$, $5< q$ and if the function $\vf:\mathbb{R} \times \R\longrightarrow \R$  belongs to the space $\M^{\frac{p}{2},\frac{q}{2}} (\mathbb{R} \times \R)$ then we have
\begin{itemize}
\item[1)] $\mathds{1}_Q\mathcal{I}_1(\vf)\in \M^{\frac{p}{\nu},\frac{q}{\nu}} (\mathbb{R} \times \R) $, with $\nu= 1-\frac{q-5}{5q}$ (remark that $0<\nu<1$).
\item[2)] $\mathds{1}_Q\mathcal{I}_1(\vf)\in \M^{\sigma, q} (\mathbb{R} \times \R) $, where $\sigma= \min\{\frac{p}{\nu},q\}$ with the same $\nu$ as before.
\end{itemize}
\end{coro}
%%%%%%%%%%%%%%%%%%%%%%%%%%%%%%%%%%%%%%%%%%%%%%%%%%%
\begin{coro}[Estimates for $\mathcal{I}_{2}$]\label{Morrey_Coro2}
Let $Q$ be a bounded set of $\mathbb{R} \times \R$ of the form given in  (\ref{SetQ}). If we have $2< p \le q$, $5<q$ and if the function $\vf:\mathbb{R} \times \R\longrightarrow \R$  belongs to the space $\M^{\frac{p}{2},\frac{q}{2}} (\mathbb{R} \times \R)$ then we have
$$\mathds{1}_Q\mathcal{I}_2(\mathds{1}_Q\vf)\in \M^{\sigma,q} (\mathbb{R} \times \R), $$ where $\sigma= \min\{\frac{p}{\nu},q\}$ with $\nu=1-\frac{q-5}{5q}$.
\end{coro}
%%%%%%%%%%%%%%%%%%%%%%%%%%%%%%%%%%%%%%%%%%%%%%%%%%%
\mysection{Gain of integrability for the variable $\vu$}\label{Secc_GainU}
Let $Q_0, Q_1$ be two parabolic balls of the form (\ref{SetQ}) such that we have the inclusions
$$Q_1\subset Q_0\subset Q.$$ 
We will start showing a small gain of integrability stated in terms of Morrey spaces. Indeed, by the hypothesis (\ref{Hypothese1}) we have
$\mathds{1}_Q\vu \in \M^{p_0,q_0}(\mathbb{R}\times\R)$ with $2<p_0\le q_0$ and $5<q_0\le 6$ and thus we can define a technical parameter $0<\nu_0<1$ such that
\begin{equation}\label{Definition_Nu0}
\nu_0=1-\frac{q_0-5}{5q_0}, 
\end{equation}
and we will prove that we have 
\begin{equation}\label{FirstMorreyIter}
\mathds{1}_{Q_1}\vu \in \M^{\sigma_0, q_0}(\mathbb{R}\times\R),
\end{equation}
where 
\begin{equation}\label{Definition_Sigma0}
\sigma_0=\min\{\frac{p_0}{\nu_0}, q_0\}.
\end{equation}
Note that since $\sigma_0>p_0$ we will obtain a (small) gain of integrability with respect of the first index of Morrey spaces and we will see later on how to iterate this process in order to deduce the first point of Theorem \ref{Theo1}.\\

We introduce now two test functions $\psi, \phi: \mathbb{R}\times \R \longrightarrow \mathbb{R}$ that belong to the space $C_0^\infty(\mathbb{R}\times \R )$ and such that 
\begin{equation}\label{Definition_FuncLocalizantes}
\begin{split}
&\psi \equiv 1\;\; on \quad Q_0\quad and \quad supp(\psi)\subset Q,\\[2mm]
&\phi \equiv 1\;\; on \;\; Q_1 \quad and \quad supp(\phi)\subset Q_0\subset Q.
\end{split}
\end{equation}
It is clear that we have $\phi(0,\cdot)=\psi(0,\cdot)=0 $ since the sets $Q_0, Q_1$ are subsets of the parabolic ball $Q$ given in (\ref{SetQ}) and $(0,0)\notin Q$, note moreover that we have the identity 
$\psi\phi =\phi$ in the whole space $\mathbb{R}\times \mathbb{R}^3$. With the help of these technical functions we can consider the function 
\begin{equation}\label{Defintion_U}
\vU= \phi\vu,
\end{equation}
hence, verifying the condition \eqref{FirstMorreyIter} amounts to prove that  $\vU \in \M^{\sigma_0, q_0}(\mathbb{R}\times \mathbb{R}^3)$ as we are only interested in the behavior of $\vu$ inside the parabolic ball $Q_1$. Now, using the localizing properties of the functions $\phi$ and $\psi$, we have in $\mathcal{S}'$ the identity
\begin{equation}\label{Definition_Ulocal}
\vU= \psi \left(\frac{1}{\Delta} \Delta(\phi \vu)\right),
\end{equation}
which can be rewritten as 
\begin{eqnarray}
\vU&=&\psi \left(\frac{1}{\Delta}\left(\phi \Delta \vu -(\Delta \phi)\vu + 2\sum_{i=1}^3 \partial_i\big((\partial _i \phi ) \vu\big)\right)\right)\notag\\
&=&\psi \left(\frac{1}{\Delta}\left(\phi \Delta \vu\right)\right) -\psi \left(\frac{1}{\Delta}\big((\Delta \phi)\vu\big)\right) + 2\sum_{i=1}^3 \psi \left(\frac{1}{\Delta}\partial_i\big((\partial _i \phi ) \vu\big)\right)\notag\\
&:=&\vU_1-\vU_2+\vU_3.\label{Def_U1U2U3}
\end{eqnarray}
Then, we shall prove that each term $\vU_1,\vU_2$ and $\vU_3$ of the expression above belong to the Morrey space $\M^{\sigma_0, q_0}$ and we will deduce from this that we have (\ref{FirstMorreyIter}). We start with the terms $\vU_2$ and $\vU_3$ which are easier to study than $\vU_1$:
%%%%%%%%%%%%%%%%%%%%%%%%%%%%%%%%%%%%%%%%%%%%%%%%%%%%%%%%%%%%%%%%
\begin{propo}\label{normL6}
Under the hypotheses of Theorem \ref{Theo1}, we have
\begin{equation}\label{NonLaplacianFirstMorrey}
\vU_2, \vU_3 \in \M^{\sigma_0, q_0}(\mathbb{R}\times\R),
\end{equation}
where $\sigma_0=\min\{\frac{p_0}{\nu_0}, q_0\}$ and $\nu_0=1-\frac{q_0-5}{5q_0}$ with $2<p_0\le q_0,\; 5<q_0\le 6$.
\end{propo}
%%%%%%%%%%%%%%%%%%%%%%%%%%%%%%%%%%%%%%%%%%%%%%%%%%%%%%%%%%%%%%%%
\noindent {\bf Proof}. We claim first that $
\vU_2, \vU_3 \in  L^\infty([0,+\infty[, L^6(\R))$. Indeed, since $\vu \in L_t^\infty L_x^2 \cap L_t^2 \dot H_x ^1(Q)$, it is easy to see that $(\Delta \phi)\vu \in L^\infty([0,+\infty[, L^{\frac{6}{5}}(\R))$ as by the H\"older inequalities we can write
\begin{eqnarray*}
\|(\Delta \phi)\vu\|_{ L_t^\infty L_x ^{\frac{6}{5}}}&\leq &\|\Delta \phi\|_{ L_t^\infty L_x ^3}\|\mathds{1}_{supp(\phi)}\vu\|_{ L_t^\infty L_x ^2}\\
&\leq &C\|\mathds{1}_{Q}\vu\|_{ L_t^\infty L_x ^2}<+\infty.
\end{eqnarray*}
Now, the embedding $L^\frac{6}{5}(\R)\subset \dot H^{-1}(\R)$ implies that
\begin{equation}\label{Info1U2}
(\Delta \phi)\vu\in L^\infty([0,+\infty[, \dot H^{-1}(\R)).
\end{equation}
Moreover, as for any $1\le i \le 3$ we have $ ( \partial _i \phi ) \vu \in L^\infty([0,+\infty[, L^2(\R))$, we obtain
\begin{equation}\label{Info1U3}
\sum_{i=1}^3\partial_i(( \partial _i \phi ) \vu)\in L^\infty([0,+\infty[, \dot H^{-1}(\R)).
\end{equation}
Thus, with the informations (\ref{Info1U2}) and (\ref{Info1U3}) above we have 
\begin{eqnarray*}
&&\vU_2=\psi \left(\frac{1}{\Delta}\left((\Delta \phi)\vu\right)\right)\in L^\infty([0,+\infty[, \dot H^{1}(\R))\quad \mbox{and}\\
&&\vU_3=2\sum_{i=1}^3  \psi \left(\frac{1}{\Delta}\partial_i \left((
 \partial _i \phi ) \vu\right)\right)\in L^\infty([0,+\infty[, \dot H^{1}(\R)),
\end{eqnarray*}
and using the classical Sobolev embedding $\dot H^{1}(\R)\subset L^{6}(\R)$ we have
$$\vU_2, \vU_3\in L^\infty([0,+\infty[, \dot H^{1}(\R))\subset L^\infty([0,+\infty[, L^{6}(\R)).$$
We remark now that due to the properties of the auxiliar function $\psi$ we have $supp(\vU_2), supp(\vU_3)\subset Q$ and since $Q$ is a bounded (parabolic) ball, we obtain that $\vU_2, \vU_3\in  L_{t,x}^6{(\mathbb{R}\times \R)}= \M^{6,6}(\mathbb{R}\times \R)$.  Finally due to the fact that $q_0\le 6$, we can conclude \eqref{NonLaplacianFirstMorrey} by using the first inequality of Lemma \ref{lemma_locindi} as we have for $i=2,3$:
$$\|\vU_i\|_{\M^{\sigma_0, q_0}}=\|\mathds{1}_Q\vU_i\|_{\M^{\sigma_0, q_0}}\leq \|\mathds{1}_Q\vU_i\|_{\M^{6, 6}}\leq\|\vU_i\|_{L^6_tL^6_x}<+\infty.$$
The proof of Proposition \ref{normL6} is finished.  \hfill $\blacksquare$
\begin{rema}\label{Remarque_GainU2U3}
Note that for the terms $\vU_2$ and $\vU_3$ we easily obtain the best control stated the expression (\ref{Conclusion1}) of the Theorem \ref{Theo1} as we have $\vU_2, \vU_3\in  L_{t,x}^6(\mathbb{R}\times \R)$. 
\end{rema}
The most technical part of the article deals with the remaining term $\vU_1$  and for this we have the following proposition
\begin{propo}\label{Proposition_32}
Under the hypotheses of Theorem \ref{Theo1}, we have
%%%%%%%%%%%%%%%%%%%%%%%%%%%%%%%%%%%%%%%%%%%%%%%%%%%
\begin{equation}\label{NonLaplacianFirstMorreyU1}
\vU_1 \in \M^{\sigma_0, q_0}(\mathbb{R}\times\R),
\end{equation}
where $\vU_1= \left(\frac{1}{\Delta}\left(\phi \Delta \vu\right)\right)$ is given in (\ref{Def_U1U2U3}) and where $\sigma_0=\min\{\frac{p_0}{\nu_0}, q_0\}$ and $\nu_0=1-\frac{q_0-5}{5q_0}$ is a small parameter with $2<p_0\le q_0,\; 5<q_0\le 6$.
\end{propo}
%%%%%%%%%%%%%%%%%%%%%%%%%%%%%%%%%%%%%%%%%%%%%%%%%%%
\noindent{\bf Proof.} In order to proof this proposition we need to study the evolution of the variable $\vU_1$ which is intimately associated to the evolution of the variable $\vu$. We note that in the first equation of \eqref{MicropolarEquations}, which is satisfied by $\vu$, the pressure $p$ is assumed to be a general object (a distribution) and to get rid of this variable $p$ in a local framework it is usual to apply the curl operator ``$\rot$''. This fact motivates the use of the following expression for $\vU_1$
\begin{equation}\label{DefinitionRotU_1}
\vU_1= \psi \left(\frac{1}{\Delta} (\phi  \Delta\vu)\right)=-\psi \left(\frac{1}{\Delta}\left(\phi (\rot [\psi \rot \vu])\right)\right),
\end{equation}
where we used the property $div(\vu)=0$ and the properties of the auxiliary functions $\phi$ and $\psi$ defined in (\ref{Definition_FuncLocalizantes}). See Lemma \ref{VectoIdentityRota} in the Appendix below for a proof of this identity.\\

It is worth noting here that in the expression (\ref{DefinitionRotU_1}) it will be enough to start our study with the quantity 
\begin{equation}\label{DefinitionRotMathcalU}
\vMU:=\rot [\psi \rot \vu],
\end{equation}
indeed, if we obtain some integrability information over $\vMU$, then it would be easy to deduce the wished Morrey information for $\vU_1$ as we have the formula
\begin{equation}\label{DefinitionRotMathcalU1}
\vU_1=\psi \left(\frac{1}{\Delta}\left(\phi \vMU\right)\right).
\end{equation}
Having this remark in mind, we will deduce the equation satisfied by $\vMU$ given in the expression (\ref{DefinitionRotMathcalU}). Thus, applying the operator ``$\vn \wedge $'' to the first equation of (\ref{MicropolarEquations}) and since we have the identity $\rot \vn p\equiv 0$ we obtain 
$$\partial_t (\vn \wedge \vu) = \Delta (\vn \wedge \vu)-\rot(\vu \cdot \vn)\vu+\tfrac{1}{2}\rot( \rot \vw)+ \rot \vf+\rot(\va\cdot \vn)\vu+\rot(\vu\cdot\vn)\va.$$
We introduce now the localizing function $\psi$ and from the previous equation we have 
\begin{eqnarray*}
\partial_t [\psi \rot \vu]&=&(\partial_t \psi)\rot \vu+ \psi\,\partial_t(\rot \vu) \\
&=&(\partial_t \psi)\rot \vu+\psi \big[\Delta (\vn \wedge \vu)-\psi\rot(\vu \cdot \vn)\vu+\tfrac{\psi}{2}\rot( \rot \vw)+ \psi\rot \vf\\
&&+\psi\rot(\va\cdot \vn)\vu+\psi\rot(\vu\cdot\vn)\va\big].
\end{eqnarray*}
Noting that we have the identity 
\begin{equation}\label{Identite_Laplacien}
\psi \Delta (\vn \wedge \vu)=\Delta (\psi\vn \wedge \vu)+(\Delta \psi) \rot \vu-2\sum_{j=1}^3\partial_j\left((\partial_j \psi)  \vn \wedge \vu\right),
\end{equation}
we can rewrite the equation above as
\begin{eqnarray*}
\partial_t [\psi \rot \vu]&=&\Delta (\psi\vn \wedge \vu)+(\partial_t\psi+\Delta \psi) (\rot \vu)-2\sum_{j=1}^3\partial_j\left((\partial_j \psi)  \vn \wedge \vu\right)\\
&&-\psi\rot(\vu \cdot \vn)\vu+\tfrac{\psi}{2}\rot( \rot \vw)+ \psi\rot \vf+\psi\rot(\va\cdot \vn)\vu+\psi\rot(\vu\cdot\vn)\va,
\end{eqnarray*}
and applying again the operator ``$\rot$'' we deduce
\begin{eqnarray*}
\partial_t (\rot[\psi \rot \vu])&=&\Delta (\rot[\psi\vn \wedge \vu])+\rot\left[(\partial_t\psi+\Delta \psi) (\rot \vu)\right]-2\rot\left[\sum_{j=1}^3\partial_j\left((\partial_j \psi)  \vn \wedge \vu\right)\right] \\
&&-\rot\left[\psi\rot(\vu \cdot \vn)\vu\right]+\rot\left[\tfrac{\psi}{2}\rot( \rot \vw)+ \psi\rot \vf\right]\\
&&+\rot\left[\psi\rot(\va\cdot \vn)\vu+\psi\rot(\vu\cdot\vn)\va\right],
\end{eqnarray*}
and since by (\ref{DefinitionRotMathcalU}) we have $\vMU=\rot [\psi \rot \vu]$:
\begin{eqnarray*}
\partial_t \vMU&=&\Delta \vMU+\rot\left[(\partial_t\psi+\Delta \psi) (\rot \vu)\right]-2\rot\left[\sum_{j=1}^3\partial_j\left((\partial_j \psi)  \vn \wedge \vu\right)\right] -\rot\left[\psi\rot(\vu \cdot \vn)\vu\right]\\
&&+\rot\left[\tfrac{\psi}{2}\rot( \rot \vw)+ \psi\rot \vf\right]+\rot\left[\psi\rot(\va\cdot \vn)\vu+\psi\rot(\vu\cdot\vn)\va\right],
\end{eqnarray*}
and this equation can be rewritten as 
\begin{equation}\label{EquaEvolU}
\partial_t \vec {\mathcal{U}}= \Delta \vec {\mathcal{U}}+\rot \vR,
\end{equation}
where the vector $\vR$ is given by the terms
\begin{eqnarray*}
\vR&=&(\partial_t\psi+\Delta \psi) (\rot \vu)-2\sum_{j=1}^3\partial_j\left((\partial_j \psi)  \vn \wedge \vu\right)+\tfrac{\psi}{2}\rot( \rot \vw)+ \psi\rot \vf\\
&& -\psi\rot(\vu \cdot \vn)\vu+\psi\rot(\va\cdot \vn)\vu+\psi\rot(\vu\cdot\vn)\va.
\end{eqnarray*}
For our purposes we need to give another expression of the three last terms above and we will use the following generic formula which is valid for $\vb, \vc$ two divergence free vector fields:
$$\psi\rot(\vb \cdot \vn)\vc=\rot\sum_{j=1}^3\partial_j(\psi b_j\vc)-\rot \sum_{j=1}^3(\partial_j\psi)b_j\vc-\sum_{j=1}^3\partial_j(\vn\psi\wedge (b_j\vc))+\sum_{j=1}^3(\partial_j\vn\psi)\wedge (b_j\vc),$$
see the Lemma \ref{Lema_VectorialIdentite1} in the Appendix below for a proof of this vectorial identity. Thus, the vector $\vR$ can be written in the following manner:
\begin{eqnarray}
\vR&=&(\partial_t\psi+\Delta \psi) (\rot \vu)-2\sum_{j=1}^3\partial_j\left((\partial_j \psi)  \vn \wedge \vu\right)+\tfrac{\psi}{2}\rot( \rot \vw)+ \psi\rot \vf\notag\\
&-&\rot\sum_{j=1}^3\partial_j(\psi u_j\vu)+\rot \sum_{j=1}^3(\partial_j\psi)u_j\vu+\sum_{j=1}^3\partial_j(\vn\psi\wedge (u_j\vu))-\sum_{j=1}^3(\partial_j\vn\psi)\wedge (u_j\vu)\label{Definition_VecteurR}\\
&+&\rot\sum_{j=1}^3\partial_j(\psi a_j\vu)-\rot \sum_{j=1}^3(\partial_j\psi)a_j\vu-\sum_{j=1}^3\partial_j(\vn\psi\wedge (a_j\vu))+\sum_{j=1}^3(\partial_j\vn\psi)\wedge (a_j\vu)\notag\\
&+&\rot\sum_{j=1}^3\partial_j(\psi u_j\va)-\rot \sum_{j=1}^3(\partial_j\psi)u_j\va-\sum_{j=1}^3\partial_j(\vn\psi\wedge (u_j\va))+\sum_{j=1}^3(\partial_j\vn\psi)\wedge (u_j\va).\notag
\end{eqnarray}
Note at this point that, by the properties of the localizing function $\psi$ defined in (\ref{Definition_FuncLocalizantes}), we have 
$$\vMU(0,x)=\rot[\psi(0,x) \rot \vu(0,x)]=0,$$ 
and thus with the equation (\ref{EquaEvolU}) we can consider the following problem
\begin{align*}
\begin{cases}
\partial_t \vec {\mathcal{U}}= \Delta \vec {\mathcal{U}}+\rot \vR,\\[3mm]
\vec {\mathcal{U}}(0,x)= 0,
\end{cases}
\end{align*}
now, by the Duhamel formula we can write
\begin{equation}\label{DuhamelFormuleU}
\vec {\mathcal{U} }(t,x)= \int _0^{t} e^{(t-s)\Delta }(\rot \vR)(s,x )ds=\rot \int _0^{t} e^{(t-s)\Delta }\vR(s,x)ds.
\end{equation}
With the help of this integral expression above for $\vMU$ and using the formula (\ref{DefinitionRotMathcalU1}) we will deduce the wished Morrey information. Indeed we have 
$$\vU_1=\psi \left(\frac{1}{\Delta}\left(\phi \vMU\right)\right)=\psi \left(\frac{1}{\Delta}\left(\phi \rot \int _0^{t} e^{(t-s)\Delta }\vR(s,x)ds\right)\right),$$
and using the identity (\ref{Definition_VecteurR}) for $\vR$ we have the following lengthy expression:
{\small
\begin{eqnarray}
\vU_1&=&\underbrace{\psi \frac{1}{\Delta}\left(\phi \rot \int _0^{t} e^{(t-s)\Delta } (\partial_t\psi+\Delta \psi) (\rot \vu)ds\right)}_{(1)}-2\underbrace{\psi \frac{1}{\Delta}\left(\phi \rot \int _0^{t} e^{(t-s)\Delta }\sum_{j=1}^3\partial_j\left((\partial_j \psi)  \vn \wedge \vu\right)ds\right)}_{(2)}\notag\\
&+&\underbrace{\psi \frac{1}{\Delta}\left(\phi \rot \int _0^{t} e^{(t-s)\Delta }\tfrac{\psi}{2}\rot( \rot \vw)ds\right)}_{(3)}+ \underbrace{\psi \frac{1}{\Delta}\left(\phi \rot \int _0^{t} e^{(t-s)\Delta }\psi\rot \vf ds\right)}_{(4)}\label{DefinitionFormuleLonguevU_1}\\
&-&\underbrace{\psi \frac{1}{\Delta}\left(\phi \rot \int _0^{t} e^{(t-s)\Delta }\rot\sum_{j=1}^3\partial_j(\psi u_j\vu)ds\right)}_{(5)}+\underbrace{\psi \frac{1}{\Delta}\left(\phi \rot \int _0^{t} e^{(t-s)\Delta} \rot \sum_{j=1}^3(\partial_j\psi)u_j\vu ds\right)}_{(6)}\notag\\
&+&\underbrace{\psi \frac{1}{\Delta}\left(\phi \rot \int _0^{t} e^{(t-s)\Delta }\sum_{j=1}^3\partial_j(\vn\psi\wedge (u_j\vu))ds\right)}_{(7)}-\underbrace{\psi \frac{1}{\Delta}\left(\phi \rot \int _0^{t} e^{(t-s)\Delta }\sum_{j=1}^3(\partial_j\vn\psi)\wedge (u_j\vu)ds\right)}_{(8)}\notag\\
&+&\underbrace{\psi \frac{1}{\Delta}\left(\phi \rot \int _0^{t} e^{(t-s)\Delta }\rot\sum_{j=1}^3\partial_j(\psi a_j\vu)ds\right)}_{(9)}-\underbrace{\psi \frac{1}{\Delta}\left(\phi \rot \int _0^{t} e^{(t-s)\Delta }\rot \sum_{j=1}^3(\partial_j\psi)a_j\vu ds\right)}_{(10)}\notag\\
&-&\underbrace{\psi \frac{1}{\Delta}\left(\phi \rot \int _0^{t} e^{(t-s)\Delta }\sum_{j=1}^3\partial_j(\vn\psi\wedge (a_j\vu))ds\right)}_{(11)}+\underbrace{\psi \frac{1}{\Delta}\left(\phi \rot \int _0^{t} e^{(t-s)\Delta }\sum_{j=1}^3(\partial_j\vn\psi)\wedge (a_j\vu)ds\right)}_{(12)}\notag\\
&+&\underbrace{\psi \frac{1}{\Delta}\left(\phi \rot \int _0^{t} e^{(t-s)\Delta }\rot\sum_{j=1}^3\partial_j(\psi u_j\va)ds\right)}_{(13)}-\underbrace{\psi \frac{1}{\Delta}\left(\phi \rot \int _0^{t} e^{(t-s)\Delta }\rot \sum_{j=1}^3(\partial_j\psi)u_j\va ds\right)}_{(14)}\notag\\
&-&\underbrace{\psi \frac{1}{\Delta}\left(\phi \rot \int _0^{t} e^{(t-s)\Delta }\sum_{j=1}^3\partial_j(\vn\psi\wedge (u_j\va))ds\right)}_{(15)}+\underbrace{\psi \frac{1}{\Delta}\left(\phi \rot \int _0^{t} e^{(t-s)\Delta }\sum_{j=1}^3(\partial_j\vn\psi)\wedge (u_j\va)ds\right).}_{(16)}\notag
\end{eqnarray}
}
With this expression for $\vU_1$, we only need to prove that each one of the previous terms belong to the Morrey space $\M^{\sigma_0, q_0}(\mathbb{R}\times\R)$ and in our study we will gather the terms that have similar structure. 
\begin{itemize}
\item We remark first that since $\vu, \vf \in L^2_t\dot{H}^1_x(Q)$ then we have $\rot \vu, \rot \vf\in L^2_tL^2_x(Q)$ and thus the terms (1) and (4) in equation (\ref{DefinitionFormuleLonguevU_1}) share the same structure and will be treated with the help of the following lemma:
\begin{lem}\label{Lemme_aux1}
Consider the parabolic ball $Q$ given in (\ref{SetQ}) and let $\vA$ be a function in $ L^2_tL^2_x(Q)$. Assume that $\psi, \phi$ are the localizing function given in (\ref{Definition_FuncLocalizantes}) and consider $\Phi$ a smooth function supported over the set $Q$. Then for $1<p\le q \le 6$ we have
\begin{equation*}
\psi \frac{1}{\Delta}\left(\phi \rot \int _0^{t} e^{(t-s)\Delta }\Phi \vA(s,\cdot )ds\right)\in \M^{p, q}(\mathbb{R}\times\R).
\end{equation*}
\end{lem}
{\bf Proof.} For the sake of simplicity let us denote $\vec{\mathcal{A}}=\displaystyle{\int _0^{t} e^{(t-s)\Delta }\Phi \vA ds}$. Since we have $1<p\leq q\le 6$ and since the function $\psi$ is supported in the parabolic ball $Q$, by the localization property given in Lemma \ref{lemma_locindi} and by the space identification $\mathcal{M}_{t,x}^{p,p}= L_{t,x}^p$, we can write
\begin{eqnarray}
\left\|\psi \frac{1}{\Delta}\left(\phi \rot \int _0^{t} e^{(t-s)\Delta }\Phi \vA(s,\cdot )ds\right)\right\|_{\M^{p, q}}=\left\|\psi \frac{1}{\Delta}\left(\phi \left(\rot \vec{\mathcal{A}}\right)\right)\right\|_{\M^{p, q}}\notag\\
\leq C\left\|\psi \frac{1}{\Delta}\left(\phi \left(\rot \vec{\mathcal{A}}\right)\right)\right\|_{\M^{6, 6}}\leq C\left\|\psi \frac{1}{\Delta}\left(\phi \left(\rot \vec{\mathcal{A}}\right)\right)\right\|_{L^6_tL^6_x}\label{EstimateL6L6Utile}\\
\leq  C \left\|\psi\right\|_{L_{t,x}^\infty}\left\|\frac{1}{\Delta}\left(\phi \left(\rot\vec{\mathcal{A}}\right)\right)\right\|_{L_{t,x}^6}.\notag
\end{eqnarray}
Now, with the embedding $\dot H^1(\R)\subset L^6(\R)$ and the properties of the negative powers of the Laplacian we obtain
$$\leq C \left\|\frac{1}{\Delta}\left(\phi \left(\rot \vec{\mathcal{A}}\right)\right)\right\|_{L_t^6\dot H _x^1}\le C \left\|\phi \left(\rot \vec{\mathcal{A}}\right)\right\|_{L_t^6\dot H _x^{-1}}\leq  C\left\|\phi \left(\rot \vec{\mathcal{A}}\right)\right\|_{L_t^6L^{\frac{6}{5}}_x},$$
where we used the dual embedding $L^{\frac{6}{5}}(\R)\subset \dot H^{-1}(\R)$ in the last estimate above. Using the fact that the function $\phi$ is supported in the (fixed and bounded) parabolic ball $Q$ we obtain by the Hölder inequality in the space variable the estimate:\\
$$\left\|\phi \left(\rot \vec{\mathcal{A}}\right)\right\|_{L_t^6L^{\frac{6}{5}}_x}\leq C\left\|\phi \left(\rot \vec{\mathcal{A}}\right)\right\|_{L_t^\infty L^{\frac{6}{5}}_x}\leq C\|\phi\|_{L_t^\infty L^3_x}\|\rot \vec{\mathcal{A}}\|_{L_t^\infty L^{2}_x}.$$
By the definition of $\vec{\mathcal{A}}$ we only have to study the quantity
\begin{eqnarray*}
\|\rot \vec{\mathcal{A}}\|_{L_t^\infty L^{2}_x}&=& \sup_{t>0}\left\|\rot\int _0^{t} e^{(t-s)\Delta } \Phi  \vA(s,\cdot ))ds\right\|_{L^2}\leq C \sup_{t>0}\left\|\int _0^{t} e^{(t-s)\Delta } \Phi  \vA(s,\cdot ))ds\right\|_{\dot{H}^1}\\
&\leq & C\|\Phi  \vA\|_{L^2_tL^2_x},
\end{eqnarray*}
where in the last line we used the properties of the heat kernel (see Lemma 7.2 of \cite{PGLR1}). To finish, it is enough to recall that the function $\Phi$ is supported in the parabolic ball $Q$ and that we have 
$\|\Phi  \vA\|_{L^2_tL^2_x}\leq C\|\vA\|_{L^2_tL^2_x(Q)}<+\infty$. We have thus proved that the quantity $\displaystyle{\psi \frac{1}{\Delta}\left(\phi \rot \int _0^{t} e^{(t-s)\Delta }\Phi \vA(s,\cdot )ds\right)}$ belongs to the space $\M^{p, q}(\mathbb{R}\times\R)$ and the proof of Lemma \ref{Lemme_aux1} is finished.\hfill $\blacksquare$\\

With this result at hand, to study the term (1) in (\ref{DefinitionFormuleLonguevU_1}) we only need to apply this lemma with $1<\sigma_0\leq q_0\leq 6$, $\Phi=(\partial_t+\Delta)\psi$ and $\vA=\rot \vu$ and for the term (4) we use $\Phi=\psi$ and $\vA=\rot \vf$. We have proven that 
$$(1)+(4)\in \M^{\sigma_0, q_0}(\mathbb{R}\times\R).$$
%%%%%%%%%%%%%%%%%%%%%%%%%%%%%%%%%%%%%%%%%%%%%%%%%%%
\item The term (2) in (\ref{DefinitionFormuleLonguevU_1}) will be studied with the following result
\begin{lem}\label{Lemme_aux2}Under the general hypotheses of Theorem \ref{Theo1} we have for $1<p\le q \le 6$ we have
\begin{equation}\label{Formule_Lemme_aux2}
\psi \frac{1}{\Delta}\left(\phi \rot \int _0^{t} e^{(t-s)\Delta }\sum_{j=1}^3\partial_j\left((\partial_j \psi)  \vn \wedge \vu\right)ds\right)\in \M^{p,q}(\mathbb{R}\times\R).
\end{equation}
\end{lem}	
%%%%%%%%%%%%%%%%%%%%%%%%%%%%%%%%%%%%%%%%%%%%%%%%%%%
{\bf Proof.}  We define $\vec{\mathcal{B}}=\displaystyle{\int _0^{t} e^{(t-s)\Delta }\sum_{j=1}^3\partial_j\left((\partial_j \psi)  \vn \wedge \vu\right)ds}$ and thus by the same previous ideas we can write
\begin{align*}
\left\|\psi \frac{1}{\Delta}\left(\phi \rot \int _0^{t} e^{(t-s)\Delta }\sum_{j=1}^3\partial_j\left((\partial_j \psi)  \vn \wedge \vu\right)ds\right) \right\|_{M_{t,x}^{\sigma_0,q_0}}=\left\|\psi\frac{1}{\Delta}(\phi \vn \wedge \vec{\mathcal{B}}) \right\|_{M_{t,x}^{\sigma_0,q_0}} \\
\leq C \left\|\psi\frac{1}{\Delta}(\phi \vn \wedge \vec{\mathcal{B}}) \right\|_{ {M_{t,x}^{6,6}}}  = C \left\|\psi\frac{1}{\Delta}(\phi \vn \wedge \vec{\mathcal{B}}) \right\|_{L^6_{t}L^6_x}.
\end{align*}
Let us define now $\vn \wedge \vec{\mathcal{B}}= \Delta \vec{\mathbb{B}}$, where $\displaystyle{\vec{\mathbb{B}} = -\sum_{j=1}^{3} \displaystyle\int_{0}^{t} e^{ (t-s)\Delta} \frac{1}{\Delta} \vn \wedge \partial_j \big( (\partial_{j}\psi)\rot \vu\big) (s,\cdot)\,ds}$ and we use the identity 
$\phi (\Delta \vec{\mathbb{B}})= \Delta(\phi\vec{\mathbb{B}})+(\Delta \phi)\vec{\mathbb{B}}-2\displaystyle{\sum_{i=1}^{3}}\partial_{i}((\partial_{i}\phi)\vec{\mathbb{B}})$ to obtain
\begin{eqnarray}\label{decomp}
\psi\Big(\frac{1}{\Delta}\big(\phi \vn \wedge \vec{\mathcal{B}}\big) \Big)&=&\psi\phi \vec{\mathbb{B}}+\psi \frac{1}{\Delta}\left((\Delta \phi) \vec{\mathbb{B}}\right)-2 \sum_{j=1}^{3} \psi \frac{\partial_{j}}{\Delta}\left(\left(\partial_{j} \phi\right) \vec{\mathbb{B}}\right).
\end{eqnarray}
We will now prove that each term of \eqref{decomp} belongs to the space $L^6_tL^6_x$. Indeed, for the first term above, by using the support properties of the functions $\psi$ and $\phi$ given in (\ref{Definition_FuncLocalizantes}) and by the Sobolev embedding  $\dot{H}^{1}(\R) \subset L^{6}(\R)$ and a standard heat kernel estimate, we get
\begin{eqnarray}
\|\psi\phi \vec{\mathbb{B}}\|_{L^6_{t}L^6_x}&\leq &C\|\vec{\mathbb{B}}\|_{L^\infty_{t} L^6_x} \leq C \| \vec{\mathbb{B}}\|_{L^\infty_{t} \dot H^1_x}\notag \\
&\leq &C\sum_{j=1}^{3} \left\| \int_{0}^{t} e^{ (t-s)\Delta} \left(\frac{1}{\Delta} \vn \wedge \partial_j \big( (\partial_{j}\psi)\rot \vu\big)\right) (s,\cdot)\,ds\right\|_{L^\infty_{t} \dot H^1_x}\notag\\
&\leq &C \sum_{j=1}^{3}\left\| \frac{1}{\Delta} \vn \wedge \partial_j \big( (\partial_{j}\psi)\rot\vu\big)\right\|_{L^2_{t} L^2_x}\leq C\sum_{j=1}^{3}\|(\partial_{j}\psi)\rot \vu\|_{L^2_{t} L^2_x}\notag\\ 
&\leq &C \| \vu \|_{L^2_{t} \dot H^1_x(Q)}<+\infty.\label{1term}
\end{eqnarray}
For the second term on the right-hand side of \eqref{decomp}, we use the embedding $\dot{H}^1(\R)\subset L^6(\R)$ and the embedding $L^{\frac{6}{5}}(\R)\subset\dot{H}^{-1}(\R)$, as well as the support properties of the function $\phi$ to get 
\begin{align*}
\left\|\psi \frac{1}{\Delta}\left((\Delta \phi) \vec{\mathbb{B}}\right) \right\|_{L^6_{t}L^6_x}
\leq C \left\| \frac{1}{\Delta}\left((\Delta \phi) \vec{\mathbb{B}}\right) \right\|_{L^6_{t}\dot{H}^1_x} 
&\leq C \left\|(\Delta \phi) \mathds{1}_{Q} \vec{\mathbb{B}} \right\|_{L^6_{t}\dot{H}^{-1}_x} \leq C \left\|(\Delta \phi) \mathds{1}_{Q} \vec{\mathbb{B}} \right\|_{L^6_{t}L^{\frac{6}{5}}_x}.
\end{align*}
Now, by the H\"older inequality we have (using again the support properties of the function $\phi$):
\begin{align*}
\left\|(\Delta \phi) \mathds{1}_{Q} \vec{\mathbb{B}} \right\|_{L^6_{t}L^{\frac{6}{5}}_x} 
\leq C \left\|(\Delta \phi)  \right\|_{L^\infty_t L^3_x} \left\| \mathds{1}_{Q} \vec{\mathbb{B}} \right\|_{L^6_{t}L^{2}_x} \leq C \| \vec{\mathbb{B}} \|_{L^\infty_{t} L^6_x},
\end{align*}
and from the estimates displayed in (\ref{1term}) we obtain $\left\|\psi \frac{1}{\Delta}\left((\Delta \phi) \vec{\mathbb{B}}\right) \right\|_{L^6_{t}L^6_x}\leq C \| \vu \|_{L^2_{t} \dot H^1_x(Q)} <+\infty$. 

For the last term of (\ref{decomp}), by the Sobolev embedding, the H\"older inequality and the inequalities (\ref{1term}) we can write 
\begin{align*}
\left\|\sum_{i=j}^{3} \psi \frac{\partial_{j}}{\Delta}\left(\left(\partial_{j} \phi\right) \vec{\mathbb{B}}\right)\right\|_{L^6_{t}L^6_x}&\leq C \sum_{j=1}^{3} \left\| \frac{\partial_{j}}{\Delta}\left(\left(\partial_{j} \phi\right) \vec{\mathbb{B}}\right)\right\|_{L^6_{t}\dot{H}^1_x} \leq C \sum_{j=1}^{3} \left\|\left(\partial_{j} \phi\right) \mathds{1}_{Q} \vec{\mathbb{B}}\right\|_{L^6_{t}L^2_x}\\
& \leq C \sum_{j=1}^{3} \left\|\partial_{j} \phi  \right\|_{L^\infty_t L^3_x} \left\| \mathds{1}_{Q} \vec{\mathbb{B}} \right\|_{L^6_{t}L^{6}_x} \leq C \| \vec{\mathbb{B}} \|_{L^\infty_{t} L^6_x}\leq C \| \vu \|_{L^2_{t} \dot H^1_x(Q)}<+\infty. 
\end{align*}
Thus gathering all the $L^6_tL^6_x$ estimates for (\ref{decomp}), we finally obtain $\psi\frac{1}{\Delta}(\phi (\vn \wedge \vec{\mathcal{B}})) \in M_{t,x}^{\sigma_0,q_0}(\mathbb{R}\times\R)$ which is the desired result. \hfill$\blacksquare$\\

We have proven, for $1<\sigma_0\leq q_0\leq 6$, that 
$$(2)\in \M^{\sigma_0, q_0}(\mathbb{R}\times\R).$$
%%%%%%%%%%%%%%%%%%%%%%%%%%%%%%%%%%%%%%%%%%%%%%%%%%%
\item The term (3) in (\ref{DefinitionFormuleLonguevU_1}) can be easily studied with the previous results. Indeed, we use the vectorial formula $\psi\rot(\rot \vw)=\rot(\psi \rot \vw)-(\vn\psi)\wedge (\rot \vw)$ to obtain
\begin{eqnarray}
\psi \frac{1}{\Delta}\left(\phi \rot \int _0^{t} e^{(t-s)\Delta }\psi\rot( \rot \vw)ds\right)=\psi \frac{1}{\Delta}\left(\phi \rot \int _0^{t} e^{(t-s)\Delta }\rot(\psi \rot \vw)ds\right) \notag\\
-\psi \frac{1}{\Delta}\left(\phi \rot \int _0^{t} e^{(t-s)\Delta }(\vn \psi)\wedge(\rot \vw)ds\right).\quad\label{TraitementDoubleRot}
\end{eqnarray}
Recall that by hypothesis we have $\rot \vw\in L^2_tL^2_x(Q)$, thus the first term in the right-hand side above can be treated using the Lemma \ref{Lemme_aux2} since this expression shares the same structure as the formula (\ref{Formule_Lemme_aux2}) above. The second term of the previous identity can be studied with the help of Lemma \ref{Lemme_aux1}. We finally obtain, for $1<\sigma_0\leq q_0\leq 6$, that  
$$(3)\in \M^{\sigma_0, q_0}(\mathbb{R}\times\R).\\[5mm]$$
\end{itemize}
%%%%%%%%%%%%%%%%%%%%%%%%%%%%%%%%%%%%%%%%%%%%%%%%%%%
It remains to study the terms (5)-(16) of (\ref{DefinitionFormuleLonguevU_1}). Recall that we have $\mathds{1}_Q\vu\in \M^{p_0, q_0}(\mathbb{R}\times\R)$ with $2<p_0\leq q_0$ and $5<q_0\leq 6$ but we also have $\mathds{1}_Q\va\in \M^{p_0, q_0}(\mathbb{R}\times\R)$ since by the Lemma \ref{lemma_locindi} we have $\|\mathds{1}_Q\va\|_{\M^{p_0, q_0}}\leq C\|\mathds{1}_Q\va\|_{\M^{6, 6}}=C\|\mathds{1}_Q\va\|_{L^6_tL^6_x}<+\infty$. Thus, since we have the same information over the variables $\vu$ and $\va$ we will perform a similar treatment of all the terms (5)-(16) following their inner structure.
%%%%%%%%%%%%%%%%%%%%%%%%%%%%%%%%%%%%%%%%%%%%%%%%%%%
\begin{itemize}
\item For the terms (5), (9) and (13) of (\ref{DefinitionFormuleLonguevU_1}) we will use the following general statement.
%%%%%%%%%%%%%%%%%%%%%%%%%%%%%%%%%%%%%%%%%%%%%%%%%%%
\begin{lem}\label{Lemme_aux3}
Consider the auxiliar functions $\psi$ and $\phi$ defined in (\ref{Definition_FuncLocalizantes}). If $\vA,\vB:\mathbb{R}\times \R\longrightarrow \R$ are two vector fields such that $\mathds{1}_{Q}\vA \in \M^{p_0,q_0}(\mathbb{R}\times \R)$ and $\mathds{1}_{Q}\vB \in \M^{p_0,q_0}(\mathbb{R}\times \R)$ with indexes $2<p_0\le q_0, 5< q_0\le 6$, then for $1\le j\le 3$ we have 
$$\psi \frac{1}{\Delta} \left(\phi\rot\int _0^{t} e^{(t-s)\Delta } \rot\partial_j(\psi A_j \vB)  ds \right)\in \M^{\sigma_0, q_0}(\mathbb{R}\times\R),$$
with $\sigma_0=\min\{\frac{p_0}{\nu_0},q_0\}$ where $\nu_0$ is given in (\ref{Definition_Nu0}).
\end{lem}
%%%%%%%%%%%%%%%%%%%%%%%%%%%%%%%%%%%%%%%%%%%%%%%%%%%
{\bf Proof.} We consider the function $\vec{\mathcal{C}}=\displaystyle{\int _0^{t} e^{(t-s)\Delta }\rot\partial_j(\psi A_i \vB)ds}$ and we define $\vec{\mathbb{C}}$ by the formula
\begin{equation*}
\vec{\mathbb{C}}= \int_0^t e^{(t-s)\Delta}\frac{1}{\Delta}\rot\rot \partial_j (\psi A_j \vB) ds 
\end{equation*}
Now, as $\rot \vec{\mathcal{C}}= \Delta \vec{\mathbb{C}}$ and by the classical identity $\phi \Delta \vec{\mathbb{C}}= \Delta (\phi \vec{\mathbb{C}})+(\Delta \phi)\vec{\mathbb{C}} -2\displaystyle{\sum_{j=1}^3}\partial_j((\partial_j \phi)\vec{\mathbb{C}})$, we can write 
{\small
\begin{equation}\label{Formule5913}
\left\|\psi \frac{1}{\Delta} \left(\phi\rot \vec{\mathcal{C}}  \right)\right\|_{ \M^{\sigma_0, q_0}}\le\left\|\psi \phi \vec{\mathbb{C}} \right\|_{ \M^{\sigma_0, q_0}}+\left\|\psi \frac{1}{\Delta}((\Delta \phi)\vec{\mathbb{C}})\right\|_{ \M^{\sigma_0, q_0}}+2 \sum_{j=1}^{3}  \left\|\psi \frac{\partial_j}{\Delta}((\partial_j \phi )\vec{\mathbb{C}})\right\|_{ \M^{\sigma_0, q_0}}.
\end{equation}
}
For the first term above, using the properties of the heat kernel and the definition of the parabolic Riesz potentials (\ref{DefinitionPotentielRiesz6}) we write
\begin{align*}
|\psi \phi \vec{\mathbb{C}}(t,x)|&=\left|\psi \phi  \int_0^t e^{(t-s)\Delta}\frac{1}{\Delta}\rot\rot \partial_j (\psi A_j \vB) ds  \right|\\
&\leq  |(\psi \phi)(t,x)| \left|\int_0^t (\partial_j e^{(t-s)\Delta})\frac{1}{\Delta}\rot\rot  (\psi A_j \vB) ds  \right|\\
&\leq C |(\psi \phi)(t,x)| \int_{\mathbb{R}}\int_{\mathbb{R}^3}
\frac{1}{(|t-s|^{\frac{1}{2}}+|x-y|)^{4}}\left|\frac{1}{\Delta}\rot\rot (\psi A_j \vB )\right|(s,y)dyds\\
&\le C |(\psi \phi)(t,x)|\mathcal{I}_1\left(\left|\frac{1}{\Delta}\rot\rot (\psi A_j \vB )\right|\right)(t,x),
\end{align*}
and using the support properties of the functions $\psi$ and $\phi$ we deduce 
\begin{equation}\label{Formula_Intermedia1}
\left\|\psi \phi  \vec{\mathbb{C}}\right\|_{ \M^{\sigma_0, q_0}}\le C\left\|\mathds{1}_{Q}\mathcal{I}_1\left(\left|\frac{1}{\Delta}\rot\rot (\psi A_j \vB )\right|\right)\right\|_{ \M^{\sigma_0, q_0}}.
\end{equation}
Now, we apply the second point of the Corollary \ref{Morrey_Coro1} with $\sigma_0=\min\{\frac{p_0}{\nu_0},q_0\}$ where $\nu_0$ is given in (\ref{Definition_Nu0}) and where $5<q_0\leq 6$ and we obtain
$$\left\|\mathds{1}_{Q}\mathcal{I}_1\left(\left|\frac{1}{\Delta}\rot\rot (\psi A_j \vB )\right|\right)\right\|_{ \M^{\sigma_0, q_0}}\leq C\left\|\frac{1}{\Delta}\rot\rot (\psi A_j \vB )\right\|_{ \M^{\frac{p_0}{2}, \frac{q_0}{2}}}\leq C\left\|\psi A_j \vB\right\|_{ \M^{\frac{p_0}{2}, \frac{q_0}{2}}},$$
where in the last estimate above we used the boundedness of the Riesz transforms in Morrey spaces. Using the H\"older inequalities we thus obtain (by the support properties of the function $\psi$):
\begin{equation}\label{Formula_Intermedia2}
\left\|\psi \phi  \vec{\mathbb{C}}\right\|_{ \M^{\sigma_0, q_0}}\leq C\left\|\mathds{1}_Q  \vA\right\|_{ \M^{p_0, q_0}}\left\|\mathds{1}_Q  \vB\right\|_{ \M^{p_0, q_0}}<+\infty.\\
\end{equation}
For the second and the third term of (\ref{Formule5913}) we will proceed as follows. First recall that the operators $\frac{1}{\Delta}$ and $\frac{\partial_j}{\Delta}$ are given by convolution with the kernels $\frac{1}{|x-y|}$ and $\frac{x_j-y_j}{|x-y|^3}$. Thus if we define the operators $T_1(\vec{\mathbb{C}})$ and $T_{2,j}(\vec{\mathbb{C}})$ for $1\leq j\leq 3$ by
\begin{eqnarray}
T_1(\vec{\mathbb{C}})(t,x)&=&\psi \frac{1}{\Delta}((\Delta \phi)\vec{\mathbb{C}})(t,x)=\psi(t,x)\int_{\R}\frac{1}{|x-y|}\Delta\phi(t,y)\vec{\mathbb{C}}(t,y)dy\label{OperateurStructure1}\\
T_{2,j}(\vec{\mathbb{C}})(t,x)&=&\psi \frac{\partial_j}{\Delta}((\partial_j \phi )\vec{\mathbb{C}})(t,x)=\psi(t,x)\int_{\R}\frac{x_j-y_j}{|x-y|^3}\partial_j\phi(t,y)\vec{\mathbb{C}}(t,y)dy.\label{OperateurStructure2}
\end{eqnarray}
we remark that the kernels associated to the operators $T_1$ and $T_{2,j}$ are bounded in $L^1(\mathbb{R}^3)$ due to the support properties of the functions $\psi$ and $\phi$. Indeed, for $T_{1}$, we have 
$$\psi(t,x) \int_{\mathbb{R}^3}\frac{1}{|x-y|}\Delta \phi(t,y)dy = \psi(t,x) \int_{Q} \frac{1}{|x-y|}\Delta\phi(t,y)dy \leq C,$$		
for almost all $x \in \mathbb{R}^3$ and
$$\Delta\phi(t,y) \int_{\mathbb{R}^3} \psi(t,x)\frac{1}{|x-y|}dx = \Delta\phi(t,y) \int_{Q} \psi(t,x)\frac{1}{|x-y|}dx \leq C,$$		
for almost all $y \in \mathbb{R}^3$.
Thus, by the Schur test, we obtain that $\|T_{1} \|_{L^1 \to L^1} \leq C$. Moreover, by the same ideas above we have $\|T_{2,j} \|_{L^1 \to L^1} \leq C$ for $1\leq j\leq 3$. Now, since the norm of $\M^{\sigma_0, q_0}$ is translation invariant, we deduce (taking into account the property $supp(\phi)\subset Q$):
\begin{equation}\label{EstimationMorreyOperateurStructure}
\begin{split}
\|T_{1}(\vec{\mathbb{C}})\|_{\M^{\sigma_0, q_0}}&=\|T_{1}(\mathds{1}_{Q}\vec{\mathbb{C}})\|_{\M^{\sigma_0, q_0}}\leq C\|\mathds{1}_{Q}\vec{\mathbb{C}}\|_{\M^{\sigma_0, q_0}}\\ 
\|T_{2,j}(\vec{\mathbb{C}})\|_{\M^{\sigma_0, q_0}}&=\|T_{2,j}(\mathds{1}_{Q}\vec{\mathbb{C}})\|_{\M^{\sigma_0, q_0}}\leq C\|\mathds{1}_{Q}\vec{\mathbb{C}}\|_{\M^{\sigma_0, q_0}}.
\end{split}
\end{equation}
Thus, coming back to the last terms of (\ref{Formule5913}) we can write 
$$\left\|\psi \frac{1}{\Delta}((\Delta \phi)\vec{\mathbb{C}})\right\|_{\M^{\sigma_0, q_0}}+2 \sum_{j=1}^{3}\left\|\psi \frac{\partial_j}{\Delta}((\partial_j \phi )\vec{\mathbb{C}})\right\|_{\M^{\sigma_0, q_0}}\leq C\|\mathds{1}_{Q}\vec{\mathbb{C}}\|_{\M^{\sigma_0, q_0}}.$$
We only need to apply the computations performed in (\ref{Formula_Intermedia1})-(\ref{Formula_Intermedia2}) to obtain
$$\left\|\psi \frac{1}{\Delta}((\Delta \phi)\vec{\mathbb{C}})\right\|_{\M^{\sigma_0, q_0}}+2 \sum_{j=1}^{3}\left\|\psi \frac{\partial_j}{\Delta}((\partial_j \phi )\vec{\mathbb{C}})\right\|_{\M^{\sigma_0, q_0}}\leq C\left\|\mathds{1}_Q  \vA\right\|_{ \M^{p_0, q_0}}\left\|\mathds{1}_Q  \vB\right\|_{\M^{p_0, q_0}}<+\infty,$$
which ends the proof of Lemma \ref{Lemme_aux3}. \hfill$\blacksquare$\\

It is straightforward to apply this results to the terms (5), (9) and (13) of (\ref{DefinitionFormuleLonguevU_1}) as we have by the hypothesis (\ref{Hypothese1}) the information $\mathds{1}_{Q}\vu\in \M^{p_0, q_0}(\mathbb{R}\times \R)$ with $2<p_0\le q_0,\; 5<q_0\le 6$ and we have $\mathds{1}_Q\va\in L^6_tL^6_x(\mathbb{R}\times \R)$ from which we easily deduce that $\mathds{1}_{Q}\va\in \M^{p_0, q_0}(\mathbb{R}\times \R)$ with the same indexes as above.\\

We have proven that 
$$(5)+(9)+(13)\in \M^{\sigma_0, q_0}(\mathbb{R}\times\R).$$
%%%%%%%%%%%%%%%%%%%%%%%%%%%%%%%%%%%%%%%%%%%%%%%%%%%
\item The terms (6), (10) and (14) of (\ref{DefinitionFormuleLonguevU_1}) we will be treated with the help of the next lemma:
\begin{lem}\label{Lemme_aux4}
Consider the functions $\psi$ and $\phi$ given in (\ref{Definition_FuncLocalizantes}). If $\vA,\vB:\mathbb{R}\times \R\longrightarrow \R$ are two vector fields such that $\mathds{1}_{Q}\vA \in \M^{p_0,q_0}(\mathbb{R}\times \R)$ and $\mathds{1}_{Q}\vB \in \M^{p_0,q_0}(\mathbb{R}\times \R)$ with indexes $2<p_0\le q_0, 5< q_0\le 6$, then for $1\le j\le 3$ we have 
$$\psi \frac{1}{\Delta}\left(\phi \rot \int _0^{t} e^{(t-s)\Delta} \rot (\partial_j\psi)A_j\vB ds\right)\in \M^{\sigma_0, q_0}(\mathbb{R}\times\R),$$
with $\sigma_0=\min\{\frac{p_0}{\nu_0},q_0\}$ where $\nu_0$ is given in (\ref{Definition_Nu0}).
\end{lem}
{\bf Proof.} Let us write $\vec{\mathcal{D}}_j(t,x)=\displaystyle{\int _0^{t} e^{(t-s)\Delta} \rot (\partial_j\psi)A_j\vB ds}$, using the properties of the heat kernel, we thus have 
\begin{equation}\label{EstimationPonctuelleRiesz61014}
|\vec{\mathcal{D}}_j(t,x)|\leq C\int_{\mathbb{R}}\int_{\R}\frac{1}{(|t-s|^\frac{1}{2}+|x-y|)^4} |(\partial_j\psi)A_j\vB(s,y)|dyds\leq  C\mathcal{I}_1(|(\partial_j\psi)A_j\vB|)(t,x),
\end{equation}
where $\mathcal{I}_1$ is the (parabolic) Riesz potential defined in (\ref{DefinitionPotentielRiesz6}). We observe now that we have 
\begin{eqnarray*}
\left\|\psi \frac{1}{\Delta}\left(\phi \rot \int _0^{t} e^{(t-s)\Delta} \rot (\partial_j\psi)A_j\vB ds\right)\right\|_{\M^{\sigma_0, q_0}}=\left\|\psi \frac{1}{\Delta}\left(\phi \rot\vec{\mathcal{D}}\right)\right\|_{\M^{\sigma_0, q_0}}\\
\leq \left\|\psi \frac{1}{\Delta}\rot\left(\phi \vec{\mathcal{D}}\right)\right\|_{\M^{\sigma_0, q_0}}+\left\|\psi \frac{1}{\Delta}\left(\vn\phi\wedge \vec{\mathcal{D}}\right)\right\|_{\M^{\sigma_0, q_0}}.
\end{eqnarray*}
At this point we remark that the operator ``$\psi \frac{1}{\Delta}\rot(\cdot)$'' in the first term above is of the same structure of the operator defined in (\ref{OperateurStructure2}) while the operator ``$\psi \frac{1}{\Delta}\left(\vn\phi\wedge \cdot\right)$'' given in the second term is of the same nature than the operator given in (\ref{OperateurStructure1}), thus by the same arguments that leaded us to (\ref{EstimationMorreyOperateurStructure}) we obtain (taking into account the support properties of the localizing function $\phi$): 
$$\left\|\psi \frac{1}{\Delta}\left(\phi \rot \int _0^{t} e^{(t-s)\Delta} \rot (\partial_j\psi)A_j\vB ds\right)\right\|_{\M^{\sigma_0, q_0}}\leq C\|\mathds{1}_{Q}\vec{\mathcal{D}}\|_{\M^{\sigma_0, q_0}}.$$
Now, with the help of the pointwise estimate (\ref{EstimationPonctuelleRiesz61014}) and by the support properties of the function $(\partial_j\psi)$ we have
$$\left\|\psi \frac{1}{\Delta}\left(\phi \rot \int _0^{t} e^{(t-s)\Delta} \rot (\partial_j\psi)A_j\vB ds\right)\right\|_{\M^{\sigma_0, q_0}}\leq C\|\mathds{1}_{Q}\mathcal{I}_1(|\mathds{1}_QA_j\vB|)\|_{\M^{\sigma_0, q_0}}.$$
We apply now the second point of Corollary \ref{Morrey_Coro1} with $\sigma_0=\min\{\frac{p_0}{\nu_0}, q_0\}$ and $\nu_0=1-\frac{q_0-5}{5q_0}$ to obtain
\begin{eqnarray*}
\left\|\psi \frac{1}{\Delta}\left(\phi \rot \int _0^{t} e^{(t-s)\Delta} \rot (\partial_j\psi)A_j\vB ds\right)\right\|_{\M^{\sigma_0, q_0}}&\leq &C\|\mathds{1}_QA_j\vB\|_{\M^{\frac{p_0}{2}, \frac{q_0}{2}}}\\
&\leq &C\left\|\mathds{1}_Q  \vA\right\|_{ \M^{p_0, q_0}}\left\|\mathds{1}_Q  \vB\right\|_{\M^{p_0, q_0}}<+\infty,
\end{eqnarray*}
where in the last line above we used the H\"older inequality in Morrey spaces. The proof of Lemma \ref{Lemme_aux4} is thus finished. \hfill$\blacksquare$\\

Since by Lemma \ref{lemma_locindi} we have $\|\mathds{1}_Q\va\|_{\M^{p_0, q_0}}\leq C\|\mathds{1}_Q\va\|_{L^6_tL^6_x}<+\infty$ and since we have by the hypothesis (\ref{Hypothese1}) that $\|\mathds{1}_Q\vu\|_{\M^{p_0, q_0}}<+\infty$, we can easily apply the previous lemma to the cases when $\vA=\vu, \va$ and $\vB=\vu, \va$ and we have proven that 
$$(6)+(10)+(14)\in \M^{\sigma_0, q_0}(\mathbb{R}\times\R).\\[5mm]$$
%%%%%%%%%%%%%%%%%%%%%%%%%%%%%%%%%%%%%%%%%%%%%%%%%%%
\item For the terms (7), (11) and (15) of (\ref{DefinitionFormuleLonguevU_1}) we use the following generic result.
\begin{lem}\label{Proposition_aux3}
If $\vA,\vB:\mathbb{R}\times \R\longrightarrow \R$ are two vector fields such that $\mathds{1}_{Q}\vA \in \M^{p_0,q_0}(\mathbb{R}\times \R)$ and $\mathds{1}_{Q}\vB \in \M^{p_0,q_0}(\mathbb{R}\times \R)$ with  $2<p_0\le q_0, 5< q_0\le 6$, then for $1\le j\le 3$ we have 
\begin{align*}
\psi \frac{1}{\Delta} \left(\phi\rot \int _0^{t} e^{(t-s)\Delta }\partial_j(\vn \psi\wedge (A_i\vB))ds \right)\in \M^{\sigma_0, q_0}(\mathbb{R}\times\R),
\end{align*}
with $\sigma_0=\min\{\frac{p_0}{\nu_0},q_0\}$ where $\nu_0$ is given in (\ref{Definition_Nu0}).
\end{lem}
{\bf Proof.} Consider the vector field $\displaystyle{\vec{\mathcal{E}}= \int _0^{t} e^{(t-s)\Delta }\partial_j( \vn\psi\wedge (A_i\vB))ds}$. By the decay properties of the heat kernel and by the support properties of the test function $\psi$, we have 
\begin{align*}
\left| \vec{\mathcal{E}}(t,x)\right|&\le \int_{0}^{t} \int_{\mathbb{R}^3} \left| \partial_j e^{(t-s)\Delta}[ \vn\psi\wedge (A_i\vB)(s,y)]\right|dyds\\
&\le C\int_{\mathbb{R}} \int_{\mathbb{R}^3} \frac{1}{(|t-s|^{\frac{1}{2}}+|x-y|)^4}\left| \vn\psi\wedge
(A_i\vB)(s,y)\right| \,dy \,ds.
\end{align*}
Using the definition of the Riesz potential given in (\ref{DefinitionPotentielRiesz6}) and by the properties of the test function $\psi$ defined in  (\ref{Definition_FuncLocalizantes}), we obtain the point-wise estimate:
\begin{equation} \label{proaux3.1}
|\vec{\mathcal{E}}(t,x)|\leq C\mathcal{I}_{1}(  \mathds{1}_{Q}|\vn\psi\wedge (A_i\vB)(t,x)|)\le 
C\mathcal{I}_{1}(\mathds{1}_{Q}| \vA(t,x)\otimes \vB(t,x)|).
\end{equation}
Recalling the identity $\phi\rot \vec{\mathcal{E}}= \rot\left(\phi \vec{\mathcal{E}}\right)-\vn\phi \wedge \vec{\mathcal{E}}$, we have
$$\left\|\psi \frac{1}{\Delta} \left(\phi\rot \vec{\mathcal{E}}  \right)\right\|_{ \M^{\sigma_0, q_0}}\leq  \left\|\psi \frac{1}{\Delta} \rot\left(\phi \vec{\mathcal{E}}  \right)\right\|_{ \M^{\sigma_0, q_0}} +   \left\|\psi \frac{1}{\Delta} \vn \phi \wedge \vec{\mathcal{E}}\right\|_{ \M^{\sigma_0, q_0}},$$
but since the operator $\frac{1}{\Delta}$ is given by convolution, we can proceed in the same fashion as in (\ref{OperateurStructure2}) for the first term above and as in (\ref{OperateurStructure1}) for the second term above and then, with the help of the estimate (\ref{EstimationMorreyOperateurStructure}), we can  deduce the inequalities $\left\|\psi \frac{1}{\Delta} \rot\left(\phi \vec{\mathcal{E}}  \right)\right\|_{ \M^{\sigma_0, q_0}} \leq C\left\|\mathds{1}_{{Q_1}} \vec{\mathcal{E}}\right\|_{ \M^{\sigma_0, q_0}}$ and  $\left\|\psi \frac{1}{\Delta} \vn \phi \wedge \vec{\mathcal{E}}\right\|_{ \M^{\sigma_0, q_0}}\leq C\left\|\mathds{1}_{{Q_1}} \vec{\mathcal{E}}\right\|_{ \M^{\sigma_0, q_0}}$. It follows that 
$$\left\|\psi \frac{1}{\Delta} \left(\phi\rot \vec{\mathcal{E}}  \right)\right\|_{ \M^{\sigma_0, q_0}}\leq C\left\|\mathds{1}_{{Q_1}} \vec{\mathcal{E}}\right\|_{ \M^{\sigma_0, q_0}}\leq C\left\|\mathds{1}_{{Q_1}} \mathcal{I}_{1}(\mathds{1}_{Q}| \vA\otimes \vB|)\right\|_{ \M^{\sigma_0, q_0}},$$
where in the last control we used the estimate (\ref{proaux3.1}) above. Thus, using the second point of the Corollary \ref{Morrey_Coro1} we have
\begin{eqnarray}
\left\|\psi \frac{1}{\Delta} \left(\phi\rot \vec{\mathcal{E}}  \right)\right\|_{ \M^{\sigma_0, q_0}} &\le & \left\|\mathds{1}_{{Q_1}} \mathcal{I}_{1}(\mathds{1}_{Q}| \vA\otimes \vB|)\right\|_{ \M^{\sigma_0, q_0}}\label{Formula_Intermedia12}\\
&\le &C \left\| \mathds{1}_{Q}| \vA\otimes  \vB|\right\|_{ \M^{\frac{p_0}{2},\frac{q_0}{2}}}\le C \left\| \mathds{1}_{Q}\vA\right\|_{ \M^{{p_0},{q_0}}}\left\| \mathds{1}_{Q}\vB\right\|_{ \M^{{p_0},{q_0}}}<+\infty\notag
\end{eqnarray}
where we applied the H\"older inequality in the framework of Morrey spaces and with this estimate the proof of Lemma \ref{Proposition_aux3} is finished. \hfill$\blacksquare$\\

Again, as we have $\|\mathds{1}_Q\va\|_{\M^{p_0, q_0}}\leq C\|\mathds{1}_Q\va\|_{L^6_tL^6_x}<+\infty$ and $\|\mathds{1}_Q\vu\|_{\M^{p_0, q_0}}<+\infty$, we can easily apply the previous lema to the cases when $\vA=\vu, \va$ and $\vB=\vu, \va$ and we have proven that 
$$(7)+(11)+(15)\in \M^{\sigma_0, q_0}(\mathbb{R}\times\R).\\[5mm]$$
%%%%%%%%%%%%%%%%%%%%%%%%%%%%%%%%%%%%%%%%%%%%%%%%%%%
\item For the terms (8), (12) and (16) of (\ref{DefinitionFormuleLonguevU_1}) we will use the following result:
\begin{lem}\label{Proposition_aux4}
If $\vA,\vB:\mathbb{R}\times \R\longrightarrow \R$ are two vector fields such that $\mathds{1}_{Q}\vA \in \M^{p_0,q_0}(\mathbb{R}\times \R)$ and $\mathds{1}_{Q}\vB \in \M^{p_0,q_0}(\mathbb{R}\times \R)$ with  $2<p_0\le q_0, 5< q_0\le 6$, then for $1\le j\le 3$ we have 
\begin{align*}
\psi \frac{1}{\Delta} \left(\phi\rot \int _0^{t} e^{(t-s)\Delta }(\partial_j\vn \psi )\wedge A_j \vB ds \right)\in \M^{\sigma_0, q_0}(\mathbb{R}\times\R),
\end{align*}
with $\sigma_0=\min\{\frac{p_0}{\nu_0},q_0\}$ where $\nu_0$ is given in (\ref{Definition_Nu0}).
\end{lem}
{\bf Proof.} As before and for the sake of simplicity we define $\displaystyle{\vec{\mathcal{F}}= \int _0^{t} e^{(t-s)\Delta }(\vn\partial_j \psi )\wedge A_j \vB ds}$. Using the decay properties of the heat equation and the definition of the parabolic Riesz potential $\mathcal{I}_{\mathfrak{2}}$ given in (\ref{DefinitionPotentielRiesz6}) we obtain the estimate
\begin{equation}\label{propoaux4.1}
|\vec{\mathcal{F}}(t,x)|\leq C\mathcal{I}_{2}(  \mathds{1}_{Q}|(\vn\partial_j \psi )\wedge A_j \vB (t,x) |)\le C\mathcal{I}_{2}(  \mathds{1}_{Q}(| \vec A(t,x)\otimes \vB(t,x)|)).   
 \end{equation}
Thus, with the same arguments displayed in the previous lemma we can write
$$\left\|\psi \frac{1}{\Delta} \left(\phi\rot \vec{\mathcal{F}}  \right)\right\|_{ \M^{\sigma_0, q_0}} \leq  \left\|\psi \frac{1}{\Delta} \rot\left(\phi \vec{\mathcal{F}}  \right)\right\|_{ \M^{\sigma_0, q_0}} +   \left\|\psi \frac{1}{\Delta} \vn \phi \wedge \vec{\mathcal{F}}\right\|_{ \M^{\sigma_0, q_0}} \leq C\left\|\mathds{1}_{{Q_1}} \vec{\mathcal{F}}\right\|_{ \M^{\sigma_0, q_0}},$$
and it follows that 
\begin{equation}\label{propoaux4.11}
\left\|\psi \frac{1}{\Delta} \left(\phi\rot \vec{\mathcal{F}}\right)\right\|_{ \M^{\sigma_0, q_0}} \leq C\left\|\mathds{1}_{{Q_1}}\mathcal{I}_{2}(  \mathds{1}_{Q}| \vec A\otimes  \vB|)\right\|_{ \M^{\sigma_0, q_0}}.
\end{equation}
At this point we invoque the Corollary \ref{Morrey_Coro2} to estimate the last term above and we can write 
$$\left\|\psi \frac{1}{\Delta} \left(\phi\rot \vec{\mathcal{F}}\right)\right\|_{ \M^{\sigma_0, q_0}}\leq C\left\| \mathds{1}_{Q}|\vec A\otimes  \vB|\right\|_{ \M^{\frac{p_0}{2},\frac{q_0}{2}}} \le C \left\| \mathds{1}_{Q}\vec A\right\|_{ \M^{{p_0},{q_0}}}\left\| \mathds{1}_{Q}\vB\right\|_{\M^{{p_0},{q_0}}}<+\infty,$$
and this finishes the proof of Lemma \ref{Proposition_aux4}.\hfill$\blacksquare$\\

As before, we can easily apply this lema to the cases when $\vA=\vu, \va$ and $\vB=\vu, \va$ and we have proven that 
$$(8)+(12)+(16)\in \M^{\sigma_0, q_0}(\mathbb{R}\times\R).\\[5mm]$$
\end{itemize}
{\bf End of the proof of Proposition \ref{Proposition_32}.} With all the previous lemmas, we have proven that all the terms of (\ref{DefinitionFormuleLonguevU_1}) belong to the Morrey space $\M^{{\sigma_0},{q_0}}(\mathbb{R}\times\R)$ and so does $\vU_1$: we have proven (\ref{NonLaplacianFirstMorreyU1}) and this ends the proof of Proposition \ref{Proposition_32}. \hfill $\blacksquare$\\
%%%%%%%%%%%%%%%%%%%%%%%%%%%%%%%%%%%%%%%%%%%%%%%%%%%
\begin{propo}\label{Proposition_GainIteratif}
Under the hypotheses of Theorem \ref{Theo1}, and recalling that $\vU$ is given in the formula (\ref{Def_U1U2U3}), we have that 
$$\vU \in L_{t,x}^{q_0}(\mathbb{R}\times\R),$$
with $5< q_0\le 6$.
\end{propo}
{\bf Proof.} Recall that following (\ref{Def_U1U2U3}) we have $\vU=\vU_1-\vU_2+\vU_3$. Note also that by Remark \ref{Remarque_GainU2U3} the terms $\vU_2$ and $\vU_3$ belong to the space $L_{t,x}^6(\mathbb{R}\times \R)$ and by the local properties of these vector fields we deduce that $\vU_2, \vU_3 \in L_{t,x}^{q_0}(\mathbb{R}\times\R)$ for $5< q_0\le 6$. It only remains to study the term $\vU_1$. For this, once we obtained that $\vU_1\in \M^{{\sigma_0},{q_0}}(\mathbb{R}\times\R)$  
with $\sigma_0=\min\{\frac{p_0}{\nu_0}, q_0\}$, then we can iterate the ideas of the Proposition \ref{Proposition_32} to obtain that $\vU_1\in \M^{\sigma_{0}', q_0}(\mathbb{R}\times\R)$ with $\sigma_{0}'= \min\{\frac{\sigma_0}{\nu_0}, q_0\}=\min\{\frac{p_0}{\nu_0^2}, q_0\}$. Since $\underset{n\to +\infty}{\lim}\frac{p_0}{\nu_0^n}=+\infty$, we easily obtain that $\vU_1 \in \M^{q_0, q_0}(\mathbb{R}\times\R)=L_{t,x}^{q_0}(\mathbb{R}\times\R)$ from which we deduce that $\vU \in L_{t,x}^{q_0}(\mathbb{R}\times\R)$.\hfill $\blacksquare$\\

To end this section, it is enough to recall that by the formula (\ref{Defintion_U}) we have the identity $\vU= \phi\vu$ and we have proven that $\phi\vu \in L_{t,x}^{q_0}(\mathbb{R}\times\R)$ with $5< q_0\le 6$ and this conclusion corresponds with the point (\ref{Conclusion1}) of Theorem \ref{Theo1}: we have thus obtained the wished gain of integrability for the variable $\vu$.
%%%%%%%%%%%%%%%%%%%%%%%%%%%%%%%%%%%%%%%%%%%%%%%%%%%
%%%%%%%%%%%%%%%%%%%%%%%%%%%%%%%%%%%%%%%%%%%%%%%%%%%
%%%%%%%%%%%%%%%%%%%%%%%%%%%%%%%%%%%%%%%%%%%%%%%%%%%
\mysection{Gain of integrability for the variable $\vw$}\label{Secc_GainOmega}
We have proven in the previous section that, under the hypotheses of Theorem \ref{Theo1}, we have the information $\mathds{1}_{Q_1}\vu \in L_{t,x}^{q_0}(\mathbb{R}\times\R)$ for some $5<q_0\leq 6$. We will now use this information in order to study the integrability properties of the second equation of (\ref{MicropolarEquations}), \emph{i.e.}:
$$\partial_t \vw = \Delta \vw +\vn div(\vw)-\vw -(\vu \cdot \vn)\vw+\frac{1}{2}\rot\vu.$$
Note that, from the generic information $\vw \in L_t^\infty L_x^2 \cap L_t^2 \dot H_x ^1(Q)$, by  an interpolation argument we easily obtain that
\begin{equation}\label{interpoW}
  \mathds{1}_Q  \vw \in L_{t,x}^{\frac{10}{3}}(\mathbb{R}\times \R)= \M ^{\frac{10}{3},\frac{10}{3}}(\mathbb{R}\times \R).
\end{equation}
Thus, with these (local) informations on $\vu$ and $\vw$ we will first prove that a (local) gain in integrability for $\vw$ is possible. For this, we introduce two auxiliar smooth functions $\varphi, \varpi: \mathbb{R}\times \R \longrightarrow \mathbb{R}$ such that for two parabolic balls $ \mathcal{Q}_a,  Q_2$ of the general form (\ref{SetQ}) that satisfy $Q_2\subset  \mathcal{Q}_a\subset Q_1$, we have
\begin{align}
\varphi &\equiv 1\;\; on \;\; \mathcal{Q}_a, \quad and \quad supp(\varphi)\subset Q_1\label{testW1}\\
\varpi &\equiv 1\;\; on \;\; Q_2, \quad and \quad supp(\varpi)\subset \mathcal{Q}_a.\label{testW2}
\end{align}
We recall that by construction we have $\varphi(0,\cdot)=\varpi(0,\cdot)=0 $ and we have the identity $\varphi \varpi=\varpi$ in the whole space. In the following proposition we will show how to obtain a first gain of integrability on $\vw$:
%%%%%%%%%%%%%%%%%%%%%%%%%%%%%%%%%%%%%%%%%%%%%%%%%%%
\begin{propo}\label{morreyw}
Under the hypotheses of Theorem \ref{Theo1} and considering the local framework stated above in (\ref{testW1})-(\ref{testW2}) we have
\begin{equation}\label{firstIteration}
\mathds{1}_{Q_2}\vw \in \M^{p,q}(\mathbb{R}\times\R)\quad with\;\;\frac{10}{3}<p\le q\le \frac{15}{4}.
\end{equation}
\end{propo}
%%%%%%%%%%%%%%%%%%%%%%%%%%%%%%%%%%%%%%%%%%%%%%%%%%%
\noindent {\bf Proof.} Due to the properties of the localizing function $\varpi$ given in (\ref{testW2}), if we define $\vW= \varpi \vw$ and if we prove that  $\vW \in   \M^{p,q}(\mathbb{R}\times\R)$ then we will easily deduce (\ref{firstIteration}) and we will focus our study in the function $\vW$. Thus, just as in the expressions (\ref{Definition_Ulocal})-(\ref{Def_U1U2U3}) above, by the properties of the functions $\varphi$ and $\varpi$ given in (\ref{testW1}) and (\ref{testW2}) we have the identity 
\begin{eqnarray*}
\vW= \varphi \left(\frac{1}{\Delta} \Delta(\varpi \vw)\right)&=& \varphi \left(\frac{1}{\Delta}(\varpi \Delta \vw)\right) -\varphi \left(\frac{1}{\Delta}((\Delta \varpi)\vw)\right) + 2\sum_{i=1}^3 \varphi \left(\frac{1}{\Delta}\partial_i\big((\partial _i \varpi )\vw\big)\right)\\
&=& \vW_1-\vW_2+\vW_3.
\end{eqnarray*}
Recall that we have $\vw\in L_t^\infty L_x^2 \cap L_t^2 \dot H_x ^1(Q)$ and following the same arguments given in the proof of the Lemma \ref{NonLaplacianFirstMorrey} we obtain that $\vW_2, \vW_3\in L_{t,x}^{6}(\mathbb{R}\times\R)=\M^{6,6}(\mathbb{R}\times\R)$ which is the strongest expected result (see Remark \ref{Remarque_GainU2U3} above): we thus only need to treat the term $\vW_1$. For this, using the Lemma \ref{VectoIdentityRota}, we can rewrite the vector field $\vW_1$ in the following manner:
\begin{eqnarray}
\vW_1&=&\varphi \left(\frac{1}{\Delta} \Delta(\varpi \vw)\right)=-\varphi \left(\frac{1}{\Delta}\left(\varpi (\rot [\varphi \rot \vw])\right)\right)+\varphi \left(\frac{1}{\Delta}\left( \varpi (\vn div(\vw ))\right)\right)\notag\\
&=&-\underbrace{\varphi \left(\frac{1}{\Delta}\left(\varpi (\rot [\varphi \rot \vw])\right)\right)}_{\vW_{1,a}}+\underbrace{\varphi \left(\frac{1}{\Delta}\left(\vn( \varpi div(\vw ))\right)\right)}_{\vW_{1,b}}-\underbrace{\varphi \left(\frac{1}{\Delta}\left((\vn \varpi) div(\vw )\right)\right)}_{\vW_{1,c}},\label{Def_W11W12W13}
\end{eqnarray}
where in the last line we used the vectorial identity $\varpi (\vn div(\vw ))=\vn  (\varpi div(\vw ))-(\vn  \varpi) div(\vw )$.\\

\noindent Now, we will prove that each term $\vW_{1,a}$, $\vW_{1,b}$ and $\vW_{1,c}$ belong to the Morrey space $\M^{p,q}(\mathbb{R}\times\R)$ with $\frac{10}{3}<p\le q\le \frac{15}{4}$.\\

In order to study the term $\vW_{1,a}$ of (\ref{Def_W11W12W13}) we introduce the variable 
\begin{equation}\label{W1A}
\vec{\mathcal{W}}_a=\varphi \rot \vw.
\end{equation}
Note that by the localizing properties of the function $\varphi$ we have $\vec{\mathcal{W}}_a(0, \cdot)=0$ and we want now to study the dynamic of this variable $\vec{\mathcal{W}}_a$. Thus, following the second equation of (\ref{MicropolarEquations}) and applying the curl operator (since we have the vectorial identity $\rot \vn div(\vw)=0$), we obtain:
$$\partial_t (\rot \vw)=\Delta (\rot \vw) -\rot \vw-\rot \big((\vu \cdot \vn)\vw\big)+\frac{1}{2}\rot(\rot \vu).$$
Now, introducing the localizing function $\varphi$ in the equation above we have
$$\partial_t [\varphi \rot \vw]=\partial_t\varphi (\rot \vw)+ \varphi\left(\Delta (\rot \vw) -\rot \vw-\rot \big((\vu \cdot \vn)\vw\big)+\frac{1}{2}\rot(\rot \vu)\right),$$
and using the identity (\ref{Identite_Laplacien}) for the term $\varphi\Delta (\rot \vw)$ we write
\begin{eqnarray*}
\partial_t [\varphi \rot \vw]&=&\partial_t\varphi (\rot \vw)+\Big(\Delta [\varphi\rot \vw]+(\Delta\varphi) (\rot \vw)-2\sum_{j=1}^3\partial_j\left((\partial_j\varphi)\rot\vw\right)\Big)\\
&&-\varphi\rot \vw-\varphi\rot \big((\vu \cdot \vn)\vw\big)+\frac{1}{2}\varphi\rot(\rot \vu).
\end{eqnarray*}
which, by the definition of the variable $\vec{\mathcal{W}}_a$ given in (\ref{W1A}) can be rewritten as 
\begin{eqnarray*}
\partial_t\vec{\mathcal{W}}_a&=&\Delta \vec{\mathcal{W}}_a+(\partial_t+\Delta\varphi) (\rot \vw)-2\sum_{j=1}^3\partial_j\left((\partial_j\varphi)\rot\vw\right)\\
&-&\varphi\rot \vw-\varphi\rot \big((\vu \cdot \vn)\vw\big)+\frac{1}{2}\varphi\rot(\rot \vu). 
\end{eqnarray*}
We need to modify the term $\varphi\rot \big((\vu \cdot \vn)\vw\big)$ and for this we write (using the fact that $div(\vu)=0$):
\begin{eqnarray*}
\varphi\rot(\vu \cdot \vn)\vw&=&\rot\sum_{j=1}^3\partial_j(\varphi u_j\vw)-\rot\sum_{j=1}^3(\partial_j\varphi) u_j\vw-\sum_{j=1}^{3}\partial_j(\vn \varphi \wedge u_j\vw)+\sum_{j=1}^{3}(\vn\partial_j \varphi )\wedge u_j \vw,
\end{eqnarray*}
thus, with this formula at hand we can write
\begin{eqnarray*}
\partial_t\vec{\mathcal{W}}_a&=&\Delta \vec{\mathcal{W}}_a+(\partial_t+\Delta\varphi) (\rot \vw)-2\sum_{j=1}^3\partial_j\left((\partial_j\varphi)\rot\vw\right)-\varphi\rot \vw\\
&-&\rot\sum_{j=1}^3\partial_j(\varphi u_j\vw)+\rot\sum_{j=1}^3(\partial_j\varphi) u_j\vw+\sum_{j=1}^{3}\partial_j(\vn \varphi \wedge u_j\vw)-\sum_{j=1}^{3}(\vn\partial_j \varphi )\wedge u_j \vw+\frac{1}{2}\varphi\rot(\rot \vu),
\end{eqnarray*}
and since we have by construction that $\vec{\mathcal{W}}_a(0, \cdot)=0$, by the Duhamel formula we can write
\begin{eqnarray*}
\vec{\mathcal{W}}_a&=& \int _0^{t} e^{(t-s)\Delta}\Bigg((\partial_t+\Delta\varphi) (\rot \vw)-2\sum_{j=1}^3\partial_j\big((\partial_j\varphi)\rot\vw\big)-\varphi\rot \vw\\
&-&\rot\sum_{j=1}^3\partial_j(\varphi u_j\vw)+\rot\sum_{j=1}^3(\partial_j\varphi) u_j\vw+\sum_{j=1}^{3}\partial_j(\vn \varphi \wedge u_j\vw)+\sum_{j=1}^{3}(\vn\partial_j \varphi )\wedge u_j \vw\\
&+&\frac{1}{2}\varphi\rot(\rot \vu)\Bigg)ds,
\end{eqnarray*}
and thus, using the definition of the variable $\vW_{1,a}$ given in (\ref{Def_W11W12W13}) we finally obtain the representation formula 
{\small
\begin{eqnarray}
\vW_{1,a}&=&\underbrace{\varphi \Bigg(\frac{1}{\Delta}\varpi \rot\int _0^{t} e^{(t-s)\Delta}(\partial_t+\Delta\varphi) (\rot \vw)ds\Bigg)}_{(1)}-2\sum_{j=1}^3\underbrace{\varphi \Bigg(\frac{1}{\Delta}\varpi\rot\int _0^{t} e^{(t-s)\Delta}\partial_j\left((\partial_j\varphi)\rot\vw\right)ds\Bigg)}_{(2)}\notag\\
&&-\underbrace{\varphi \Bigg(\frac{1}{\Delta}\varpi \rot\int _0^{t} e^{(t-s)\Delta}\varphi\rot \vw ds\Bigg)}_{(3)}-\sum_{j=1}^3\underbrace{\varphi \Bigg(\frac{1}{\Delta}\varpi\rot\int _0^{t} e^{(t-s)\Delta}\rot\partial_j(\varphi u_j\vw)ds\Bigg)}_{(4)}\label{FormulaW1A}\\
&&+\sum_{j=1}^3\underbrace{\varphi \Bigg(\frac{1}{\Delta}\varpi \rot\int _0^{t} e^{(t-s)\Delta}\rot(\partial_j\varphi) u_j\vw ds\Bigg)}_{(5)}+\sum_{j=1}^{3}\underbrace{\varphi \Bigg(\frac{1}{\Delta}\varpi \rot\int _0^{t} e^{(t-s)\Delta}\partial_j(\vn \varphi \wedge u_j\vw)ds\Bigg)}_{(6)}\notag\\
&&-\sum_{j=1}^{3}\underbrace{\varphi \Bigg(\frac{1}{\Delta}\varpi \rot\int _0^{t} e^{(t-s)\Delta}(\vn\partial_j \varphi )\wedge u_j \vw ds\Bigg)}_{(7)}+\frac{1}{2}\underbrace{\varphi \Bigg(\frac{1}{\Delta}\varpi \rot\int _0^{t} e^{(t-s)\Delta}\varphi\rot(\rot \vu)ds\Bigg)}_{(8)}\notag
\end{eqnarray}
}
and thus we need to prove that all the previous terms belong to the Morrey space $\M^{p,q}(\mathbb{R}\times\R)$ with $\frac{10}{3}<p\le q\le \frac{15}{4}$ and we will study them in the following points:
\begin{itemize}
\item The terms (1) and (3) of (\ref{FormulaW1A}) can be treated in the same manner since we have $\rot \vw\in L^2_tL^2_x(Q)$. Indeed, by applying the Lemma \ref{Lemme_aux1} with $p=q=6$ we easily obtain that the terms (1) and (3) belong to the space $L^6_tL^6_x(\mathbb{R}\times\R)$ (see the estimate (\ref{EstimateL6L6Utile}) above) and due to the localizing properties of the function $\varphi$ given in (\ref{testW1}) we easily deduce (by Lemma \ref{lemma_locindi}) that 
$$(1)+(3)\in\M^{p,q}(\mathbb{R}\times\R)\qquad \mbox{with} \;\tfrac{10}{3}<p\le q\le \tfrac{15}{4}.$$
\item For the term (2) of (\ref{FormulaW1A}), since $\rot \vw\in L^2_tL^2_x(Q)$ we can apply Lemma \ref{Lemme_aux2} and we easily deduce that (2) belongs to the space $L^6_tL^6_x(\mathbb{R}\times\R)$, from which we obtain 
$$(2)\in\M^{p,q}(\mathbb{R}\times\R)\qquad \mbox{with} \;\tfrac{10}{3}<p\le q\le \tfrac{15}{4}.$$

\item The term (4) of (\ref{FormulaW1A}) will be treated as follows. Following the same ideas given in the Lemma \ref{Lemme_aux3} (see in particular the estimates
(\ref{Formula_Intermedia1})-(\ref{EstimationMorreyOperateurStructure})), by the support properties of the functions $\varphi$ and $\varpi$  we can write
\begin{equation}\label{Formula41}
\left\|\varphi \frac{1}{\Delta}\varpi\rot\int _0^{t} e^{(t-s)\Delta}\rot\partial_j(\varphi u_j\vw)ds\right\|_{\M^{p,q}}\leq  C\left\|\mathds{1}_{Q_2}\mathcal{I}_1\left(\left|\frac{1}{\Delta}\rot\rot (\varphi \vu\otimes\vw)\right|\right)\right\|_{ \M^{p, q}}.
\end{equation}
 Since $\tfrac{10}{3}<p\le q\le \tfrac{15}{4}$, if we consider $2<q_1\leq \frac{15}{7}$ then we have $1<\frac{5}{q_1}$ and if we define $\nu_1=1-\frac{q_1}{5}$, we have $\tfrac{10}{3}<\frac{q_1}{\nu_1}\le \tfrac{15}{4}$. Thus, by properties given in Lemma \ref{lemma_locindi} we can write
$$\left\|\mathds{1}_{Q_2}\mathcal{I}_1\left(\left|\frac{1}{\Delta}\rot\rot (\varphi \vu\vw)\right|\right)\right\|_{ \M^{p, q}}\leq \left\|\mathds{1}_{Q_2}\mathcal{I}_1\left(\left|\frac{1}{\Delta}\rot\rot (\varphi \vu\otimes\vw)\right|\right)\right\|_{ \M^{\frac{q_1}{\nu_1}, \frac{q_1}{\nu_1}}},$$
and using the Lemma \ref{Lemme_Hed} we deduce the estimate
$$\left\|\mathds{1}_{Q_2}\mathcal{I}_1\left(\left|\frac{1}{\Delta}\rot\rot (\varphi \vu\otimes\vw)\right|\right)\right\|_{ \M^{\frac{q_1}{\nu_1}, \frac{q_1}{\nu_1}}}\leq \left\|\frac{1}{\Delta}\rot\rot (\varphi \vu\otimes\vw)\right\|_{ \M^{q_1, q_1}}\leq \|\varphi \vu\otimes\vw\|_{L^{q_1}_{t,x}},$$
where we used the boundedness of the Riesz transforms in Morrey spaces and the identification between Morrey and Lebesgue spaces. Thus, by the usual H\"older inequalities with $\frac{1}{q_1}=\frac{1}{q_0}+\frac{3}{10}$ (note that the condition $2<q_1\leq \tfrac{15}{7}$ stated above is related to the previous identity and the condition $5<q_0\leq 6$) we obtain (since $Q_1\subset Q$) 
\begin{equation}\label{Formula42}
\|\mathds{1}_{Q_1} \vu\otimes  \vw\|_{L^{q_1}_{t,x}}\leq \|\mathds{1}_{Q_1} \vu\|_{L_{t,x}^{q_0}}\|\mathds{1}_{Q} \vw\|_{L_{t,x}^{\frac{10}{3}}}<+\infty,
\end{equation}
which is a bounded quantity since we have already proven $\|\mathds{1}_{Q_1} \vu\|_{L_{t,x}^{q_0}}<+\infty$ and we have by interpolation the information $\|\mathds{1}_{Q} \vw\|_{L_{t,x}^{\frac{10}{3}}}<+\infty$. Thus, we obtain that 
$$(4)\in\M^{p,q}(\mathbb{R}\times\R)\qquad \mbox{with} \;\tfrac{10}{3}<p\le q\le \tfrac{15}{4}.$$

\item For the term (5) of (\ref{FormulaW1A}), by the same arguments displayed in the Lemma \ref{Proposition_aux3} (see the inequality (\ref{Formula_Intermedia12}) above) we obtain the estimate 
\begin{equation}\label{Formula51}
\left\|\varphi \frac{1}{\Delta}\varpi \rot\int _0^{t} e^{(t-s)\Delta}\rot(\partial_j\varphi) u_j\vw ds\right\|_{\M^{p,q}}\leq C\left\|\mathds{1}_{Q_2}\mathcal{I}_1(\mathds{1}_{Q_1}|\vu\otimes \vw|)\right\|_{\M^{p,q}},
\end{equation}
following the ideas of the previous point (see formula (\ref{Formula42})) we can write
$$\left\|\mathds{1}_{Q_2}\mathcal{I}_1(\mathds{1}_{Q_1}|\vu\otimes \vw|)\right\|_{\M^{p,q}}\leq \left\|\mathcal{I}_1(\mathds{1}_{Q_1}|\vu\otimes \vw|)\right\|_{ \M^{\frac{q_1}{\nu_1}, \frac{q_1}{\nu_1}}}\le  \|\mathds{1}_{Q_1} \vu\|_{L_{t,x}^{q_0}}\|\mathds{1}_{Q} \vw\|_{L_{t,x}^{\frac{10}{3}}}<+\infty,$$
and we have proven 
$$(5)\in\M^{p,q}(\mathbb{R}\times\R)\qquad \mbox{with} \;\tfrac{10}{3}<p\le q\le \tfrac{15}{4}.$$

\item The term (6) of (\ref{FormulaW1A}) shares the same structure of the term (5) and thus, by the same arguments we obtain
$$(6)\in\M^{p,q}(\mathbb{R}\times\R)\qquad \mbox{with} \;\tfrac{10}{3}<p\le q\le \tfrac{15}{4}.$$

\item For the term (7) of (\ref{FormulaW1A}), following the same ideas displayed in the proof of the Lemma \ref{Proposition_aux4} (see in particular the estimates (\ref{propoaux4.1})-(\ref{propoaux4.11})) we can write
\begin{equation}\label{Formula71}
\left\|\varphi \Bigg(\frac{1}{\Delta}\varpi \rot\int _0^{t} e^{(t-s)\Delta}(\vn\partial_j \varphi )\wedge u_j \vw ds\Bigg)\right\|_{ \M^{p,q}}\leq C \|\mathds{1}_{Q_2}\mathcal{I}_2(| \mathds{1}_{Q_1}\vu\otimes  \vw|)\|_{ \M^{p,q}}.
\end{equation}
Again, since $\tfrac{10}{3}<p\le q\le \tfrac{15}{4}$ and $2<q_1\leq \frac{15}{7}$, we have $2<\frac{5}{q_1}$ and if we define  $\nu_1=1-\frac{2q_1}{5}$ then we have $10<\frac{q_1}{\nu_1}\le 15$ (note that with these values we obtain directly the wished indexes of integrability). Thus by the Lemma \ref{lemma_locindi} and by Lemma \ref{Lemme_Hed} we have 
$$C \|\mathds{1}_{Q_2}\mathcal{I}_2(| \mathds{1}_{Q_1}\vu\otimes  \vw|)\|_{ \M^{p,q}}\leq C \|\mathcal{I}_2(| \mathds{1}_{Q_1}\vu\otimes  \vw|)\|_{ \M^{\frac{q_1}{\nu_1},\frac{q_1}{\nu_1}}}\leq  C\| \mathds{1}_{Q_1}\vu\otimes  \vw\|_{ \M^{q_1,q_1}}=C\|\mathds{1}_{Q_1} \vu\otimes  \vw\|_{L^{q_1}_{t,x}},$$
and once we have this estimate, recalling the inequalities (\ref{Formula42}) we can easily conclude that
$$(7)\in\M^{p,q}(\mathbb{R}\times\R)\qquad \mbox{with} \;\tfrac{10}{3}<p\le q\le \tfrac{15}{4}.$$

\item The term (8) of (\ref{FormulaW1A}) can be treat just as in the formula (\ref{TraitementDoubleRot}) above and since $\rot \vu\in L^2_tL^2_x(Q)$ and we obtain 
$$(8)\in\M^{p,q}(\mathbb{R}\times\R)\qquad \mbox{with} \;\tfrac{10}{3}<p\le q\le \tfrac{15}{4}.\\[3mm]$$
With the estimates of the terms (1)-(9) of (\ref{FormulaW1A}) we have proven that the term $\vec{W}_{1,a}$ given in (\ref{Def_W11W12W13}) belongs to the space $\M^{p,q}(\mathbb{R}\times\R)$ with $\tfrac{10}{3}<p\le q\le \tfrac{15}{4}$.\\[5mm]
\end{itemize}
%%%%%%%%%%%%%%%%%%%%%%%%%%%%%%%%%%%%%%%%%%%%%%%%%%%
\noindent We treat now the term $\vec{W}_{1,b}$ defined in (\ref{Def_W11W12W13}) and for this we define the variable 
$$\mathcal{W}_b= \varpi div(\vw),$$
which satisfies the following system (derived from the second equation of (\ref{MicropolarEquations})):
$$\partial_t\mathcal{W}_b=2 \Delta \mathcal{W}_b+(\partial_t\varpi+2\Delta\varpi-\varpi)div(\vw)-4\sum_{j=1}^3\partial_j((\partial_j \varpi)div(\vw))-\varpi div((\vu\cdot \vn)\vw),$$
where we used the identity (\ref{Identite_Laplacien}) for the term $\varpi \Delta div(\vw)$ and the fact that $div(\vu)=0$. Noting that, by the support properties of the function $\varpi$, we have $\mathcal{W}_b(0, \cdot)=0$, we can write by the Duhamel formula:
\begin{equation}\label{funcW}
\mathcal{W}_b=\int_0^t e^{2(t-s)\Delta}\Big((\partial_t\varpi+2\Delta\varpi-\varpi)div(\vw)-4\sum_{j=1}^3\partial_j((\partial_j \varpi)div(\vw))-\varpi div((\vu\cdot \vn)\vw)\Big) ds,
\end{equation}
thus we obtain the expression
\begin{eqnarray}
\vec{W}_{1,b}&=&\underbrace{\varphi \frac{1}{\Delta}\vn\int_0^t e^{2(t-s)\Delta}\Big((\partial_t\varpi+2\Delta\varpi-\varpi)div(\vw)\Big)ds}_{(1)} -4\sum_{j=1}^3\underbrace{\varphi \frac{1}{\Delta}\vn\int_0^t e^{2(t-s)\Delta}\Big(\partial_j((\partial_j \varpi)div(\vw))\Big) ds}_{(2)}\notag\\
&&-\underbrace{\varphi \frac{1}{\Delta}\vn\int_0^t e^{2(t-s)\Delta}\Big(\varpi div((\vu\cdot \vn)\vw)\Big) ds}_{(3)},\label{TermS}
\end{eqnarray}
and we will prove that each term above belongs to the Morrey space $\M^{p,q}(\mathbb{R}\times\R)$ with indexes $\tfrac{10}{3}<p\le q\le \tfrac{15}{4}$.
\begin{itemize}
\item For the first term of (\ref{TermS}) we denote by $\phi=(\partial_t\varpi+2\Delta\varpi-\varpi)$. Since we have $p\le q\le \tfrac{15}{4}<6$, by Lemma \ref{lemma_locindi} we can write
\begin{eqnarray*}
\left\|\psi \left(\frac{1}{\Delta}\left(\vn \int_0^t e^{2(t-s)\Delta}\phi )div(\vw)ds\right)\right)\right\|_{\M^{p, q}}&\le &\left\|\psi \left(\frac{1}{\Delta}\left(\vn \int_0^t e^{2(t-s)\Delta} \phi div(\vw)ds\right)\right)\right\|_{L_{t,x}^6}\\
&\leq & \|\psi\|_{L^6_tL^\infty_x} \left\|\frac{1}{\Delta}\left(\vn \int_0^t e^{2(t-s)\Delta} \phi div(\vw)ds\right)\right\|_{L^\infty_tL^6_x}, 
\end{eqnarray*}
and by the Sobolev embedding $\dot H^1(\R)\subset L^6(\R)$ we obtain 
\begin{eqnarray*}
\left\|\psi \left(\frac{1}{\Delta}\left(\vn \int_0^t e^{2(t-s)\Delta}\phi )div(\vw)ds\right)\right)\right\|_{\M^{p, q}}&\leq &C\left\|\frac{1}{\Delta}\left(\vn \int_0^t e^{2(t-s)\Delta} \phi div(\vw)ds\right)\right\|_{L^\infty_t\dot{H}^1_x}\\
&\leq & C\left\|\int_0^t e^{2(t-s)\Delta} \phi div(\vw)ds\right\|_{L_t^\infty L _x^{2}}.
\end{eqnarray*}
Then, using the properties of the heat kernel, the dual embedding $\dot H^{-1}(\R)\subset L^{\frac{6}{5}}(\R)$, the H\"older inequality and the properties of the function $\phi$, we have 
\begin{eqnarray*}
\left\|\int_0^t e^{2(t-s)\Delta} \phi div(\vw)ds\right\|_{L_t^\infty L_x^{2}}&\leq &C\|\phi div(\vw)\|_{L_t^2 \dot{H}^{-1}_x}\leq C\|\phi div(\vw)\|_{L_t^2 L^{\frac{6}{5}}_x}\\
&\leq &C\|\phi\|_{L^\infty_tL^3_x}\|\mathds{1}_Q div(\vw)\|_{L_t^2 L^{2}_x}<+\infty, 
\end{eqnarray*}
from which we deduce that the quantity (1) of (\ref{TermS})  belongs to the wished space $\M^{p,q}(\mathbb{R}\times\R)$ with $\tfrac{10}{3}<p\le q\le \tfrac{15}{4}$.

\item For the second term of (\ref{TermS}), by the same arguments as above we obtain
\begin{eqnarray*}
\left\|\varphi \frac{1}{\Delta}\vn\int_0^t e^{2(t-s)\Delta}\Big(\partial_j((\partial_j \varpi)div(\vw))\Big) ds\right\|_{\M^{p, q}}\leq C\left\|\int_0^te^{2(t-s)\Delta}\Big(\partial_j((\partial_j \varpi)div(\vw))\Big) ds\right\|_{L^\infty_tL^2_x}\\
\leq  C\left\|\int_0^te^{2(t-s)\Delta}(\partial_j \varpi)div(\vw)ds\right\|_{L^\infty_t\dot{H}^1_x}\leq \|(\partial_j \varpi)div(\vw)\|_{L^2_tL^2_x}\leq C\|\mathds{1}_Q div(\vw)\|_{L_t^2 L^{2}_x}<+\infty, 
\end{eqnarray*}
and thus  $(2) \in \M^{p,q}(\mathbb{R}\times\R)$ with $\tfrac{10}{3}<p\le q\le \tfrac{15}{4}$.
\item For the last term of (\ref{TermS}), we need to study the quantity 
\begin{equation}\label{TermSFinal}
\left\|\varphi \frac{1}{\Delta}\vn\int_0^t e^{2(t-s)\Delta}\Big(\varpi div((\vu\cdot \vn)\vw)\Big) ds\right\|_{\M^{p, q}},
\end{equation}
but we remark that since we have $div(\vu)=0$ we can write $\varpi div((\vu\cdot \vn)\vw)=\varpi div(div(\vw\otimes \vu))$ and then we should deal with terms of the form $\varpi \partial_i(\partial_j(u_k w_\ell))$ for $1\leq i,j,k,\ell\leq 3$, which can we rewritten in the following manner
\begin{equation}\label{TermSFinal1}
\varpi \partial_i(\partial_j(u_k w_\ell))=\partial_i(\partial_j(\varpi u_k w_\ell))-\partial_i((\partial_j\varpi)u_kw_\ell)-\partial_j((\partial_i\varpi)u_kw_\ell)-(\partial_i\partial_j\varpi)u_kw_\ell.
\end{equation}
Thus, to study the quantity (\ref{TermSFinal}) we first consider the expression
$$\left\|\varphi \frac{1}{\Delta}\vn\int_0^t e^{2(t-s)\Delta}\partial_i(\partial_j(\varpi u_k w_\ell)) ds\right\|_{\M^{p, q}},$$
which is of the same shape of the left-hand side of the formula (\ref{Formula41}) and thus, by the same arguments displayed in (\ref{Formula41})-(\ref{Formula42}) we easily obtain that this quantity is bounded in the wished Morrey space $\M^{p,q}(\mathbb{R}\times\R)$.  Since the second and the third term of (\ref{TermSFinal1}) are of the same structure, we only study one of them and we write 
$$\left\|\varphi \frac{1}{\Delta}\vn\int_0^t e^{2(t-s)\Delta}\partial_i((\partial_j\varpi)u_kw_\ell) ds\right\|_{\M^{p, q}}.$$
This term can be treated just as the left-hand side of (\ref{Formula51}) and we obtain that this quantity is bounded. Finally, for the last term of (\ref{TermSFinal1}) we write
$$\left\|\varphi \frac{1}{\Delta}\vn\int_0^t e^{2(t-s)\Delta}(\partial_i\partial_j\varpi)u_kw_\ell ds\right\|_{\M^{p, q}},$$
and following the computations made in (\ref{Formula71}) we also obtain that it is bounded in the corresponding Morrey space. We have studied all the terms of (\ref{TermSFinal1}) and then we deduce that the quantity  (\ref{TermSFinal}) belongs to the Morrey space $\M^{p,q}(\mathbb{R}\times\R)$ with $\tfrac{10}{3}<p\le q\le \tfrac{15}{4}$.\\[5mm]
\end{itemize}
%%%%%%%%%%%%%%%%%%%%%%%%%%%%%%%%%%%%%%%%%%%%%%%%%%%
We need now to study the last term of (\ref{Def_W11W12W13}), which is given by $\vW_{1,c}=\varphi \frac{1}{\Delta}\big((\vn \varpi) div(\vw )\big)$, and we study the quantity 
$$\left\|\varphi \frac{1}{\Delta}\left((\vn\varpi) div(\vw)\right)\right\|_{\M^{p, q}}.$$
Let us introduce the function $\vec\phi=\vn\varpi$, and by the identity $\phi_j \partial_k(w_\ell)=\partial_k(\phi_j w_\ell)-(\partial_k\phi_j)w_\ell$, it is enough to treat the quantities 
\begin{equation}\label{FormulaW13}
\left\|\varphi \frac{1}{\Delta}\left(\partial_k(\phi_j w_\ell)\right)\right\|_{\M^{p, q}}, \qquad \left\|\varphi \frac{1}{\Delta}\left((\partial_k\phi_j)w_\ell\right)\right\|_{\M^{p, q}},
\end{equation}
for all $1\leq j,k,\ell\leq 3$. For the first term above, using the support properties of the function $\varphi$, the fact that $p\leq q<6$ and the Sobolev embedding $\dot H ^1(\R)\subset L^6 (\R)$, we write
$$\left\|\varphi \frac{1}{\Delta}\left(\partial_k(\phi_j w_\ell)\right)\right\|_{\M^{p, q}}\leq C\left\| \varphi\frac{1}{\Delta}\left(\partial_k(\phi_j w_\ell)\right)\right\|_{L^6_tL^6_x}\leq C\left\|\frac{1}{\Delta}\left(\partial_k(\phi_j w_\ell)\right)\right\|_{L^\infty_t\dot{H}^1_x}\leq C\|\phi_j w_\ell\|_{L^\infty_tL^2_x}<+\infty,$$
since we have that $\|\mathds{1}_Q\vw\|_{L^\infty_tL^2_x}<+\infty$. For the second term of (\ref{FormulaW13}) we have
$$\left\|\varphi \frac{1}{\Delta}\left((\partial_k\phi_j)w_\ell\right)\right\|_{\M^{p, q}}\leq C\left\| \varphi\frac{1}{\Delta}\left((\partial_k\phi_j)w_\ell\right)\right\|_{L^6_tL^6_x}\leq C\left\|\frac{1}{\Delta}\left((\partial_k\phi_j)w_\ell\right)\right\|_{L^\infty_t\dot{H}^1_x}\leq C\|(\partial_k\phi_j)w_\ell\|_{L^\infty_t\dot{H}^{-1}_x},$$
and by the embedding  $L^{\frac{6}{5}}(\R)\subset \dot{H}^{-1}(\R)$, we can write
$$\left\|\varphi \frac{1}{\Delta}\left((\partial_k\phi_j)w_\ell\right)\right\|_{\M^{p, q}}\leq C\|(\partial_k\phi_j)w_\ell\|_{L^\infty_tL^{\frac{6}{5}}_x}\leq \|\partial_k\phi_i\|_{L_t^\infty L_x^{3}}\|\mathds{1}_Qw_\ell\|_{L_t^\infty L_x^{2}}<+\infty.$$
We have proven that all the terms of (\ref{FormulaW13}) belong to the Morrey space $\M^{p,q}(\mathbb{R}\times\R)$ with $\tfrac{10}{3}<p\le q\le \tfrac{15}{4}$ and so does the quantity $\vW_{1,c}$.\\[5mm]

We have proven so far that all the term in the expression (\ref{Def_W11W12W13}) belong to the Morrey space  $\M^{p,q}(\mathbb{R}\times\R)$ and this ends the proof of the Proposition \ref{morreyw} $\hfill \blacksquare$\\

%%%%%%%%%%%%%%%%%%%%%%%%%%%%%%%%%%%%%%%%%%%%%%%%%%%
We now return to the proof of the second point of Theorem \ref{Theo1}: we have proven so far the local information $\mathds{1}_{Q_2}\vw\in \M^{p, q}(\mathbb{R}\times\R)$ with $\frac{10}{3}<p\leq q\le \frac{15}{4}$ and in particular, due to the Lemma \ref{lemma_locindi} and the identification of Morrey and Lebesgue spaces, we have $\mathds{1}_{Q_2}\vw\in \M^{q,q}(\mathbb{R}\times\R)=L^{q}_{t,x}(\mathbb{R}\times\R)$ with $\frac{10}{3}<q\le \frac{15}{4}$ which is a gain of information with respect to (\ref{interpoW}). It is then enough to replace this information in the previous arguments (reducing if necessary the support of the auxiliary functions) to obtain a higher integrability control, and by suitable iterations we finally obtain that 
$$\mathds{1}_{Q_2}\vw \in L_{t,x}^{q_0}(\mathbb{R}\times\R),$$
for some $5<q_0\leq 6$, which corresponds with the wished estimate (\ref{Conclusion2}). The proof of Theorem \ref{Theo1} is now complete. \hfill$\blacksquare$
%%%%%%%%%%%%%%%%%%%%%%%%%%%%%%%%%%%%%%%%%%%%%%%%%%%
\section*{Appendix}
%%%%%%%%%%%%%%%%%%%%%%%%%%%%%%%%%%%%%%%%%%%%%%%%%%%
\begin{lem}\label{VectoIdentityRota}
If $\phi, \psi$ are the test functions given in (\ref{Definition_FuncLocalizantes}), and $\vu$ is a regular enough vector field, we have
 \begin{equation*}
      \psi \left(\frac{1}{\Delta} (\phi  \Delta\vu)\right)= 
      -\psi \left(\frac{1}{\Delta}\left(\phi (\rot [\psi \rot \vu])\right)\right)
      +\psi \left(\frac{1}{\Delta}\left(\phi(\vn div(\vu ))\right)\right).
 \end{equation*}
\end{lem}
%%%%%%%%%%%%%%%%%%%%%%%%%%%%%%%%%%%%%%%%%%%%%%%%%%%
\noindent{\bf Proof.}
We have
\begin{align*}
    \rot [\psi \rot \vu]&=\psi \rot (\rot \vu)+\vn \psi\wedge(\rot \vu)=\psi(\vn div(\vu )-\Delta \vu)+\vn \psi\wedge(\rot \vu).
\end{align*}
Moreover, by the support properties of $\phi$ and $\psi$ we have $\phi\vn\psi\equiv 0$ and $\psi \phi=\phi$ on the support of $\phi$. Then, the second term in the identity above disappear when we multiply by $\phi$ and we obtain
\begin{equation*}
    -\phi \rot [\psi \rot \vu]+\phi\vn div(\vu )=\phi\Delta \vu,
\end{equation*}
and from this identity it is easy to recover the wished result. \hfill $\blacksquare$
%%%%%%%%%%%%%%%%%%%%%%%%%%%%%%%%%%%%%%%%%%%%%%%%%%%
\begin{lem}\label{Lema_VectorialIdentite1}
Assume that $\vb$ and $\vc$ are two vector fields such that $div(\vb)=div(\vc)=0$ and assume that $\psi$ is the smooth function given in (\ref{Definition_FuncLocalizantes}). Then we have the following vectorial identity
$$\psi\rot(\vb \cdot \vn)\vc=\rot\sum_{j=1}^3\partial_j(\psi b_j\vc)-\rot \sum_{j=1}^3(\partial_i\psi)b_j\vc-\sum_{j=1}^3\partial_j(\vn\psi\wedge (b_j\vc))+\sum_{j=1}^3(\partial_j\vn\psi)\wedge (b_j\vc).$$
\end{lem}
%%%%%%%%%%%%%%%%%%%%%%%%%%%%%%%%%%%%%%%%%%%%%%%%%%%
\noindent{\bf Proof.} We start with the formula $\psi\rot(\vb \cdot \vn)\vc=\rot\big(\psi (\vb \cdot \vn)\vc\big)-\vn\psi\wedge \big((\vb \cdot \vn)\vc\big)$ which can be rewritten as 
$$\psi\rot(\vb \cdot \vn)\vc=\rot\sum_{j=1}^3\partial_j(\psi b_j\vc)- \rot\sum_{j=1}^3\partial_j(\psi b_j)\vc-\sum_{j=1}^3\vn\psi\wedge b_j\partial_j\vc.$$
Using now the fact that $div(\vb)=div(\vc)=0$ in the second and the third term of the right-hand side above we obtain 
$$\psi\rot(\vb \cdot \vn)\vc=\rot\sum_{j=1}^3\partial_j(\psi b_j\vc)- \rot\sum_{j=1}^3(\partial_j\psi) b_j\vc-\sum_{j=1}^3\vn\psi\wedge \partial_j(b_j\vc),$$
from which we easily deduce that 
$$\psi\rot(\vb \cdot \vn)\vc=\rot\sum_{j=1}^3\partial_j(\psi b_j\vc)-\rot \sum_{j=1}^3(\partial_i\psi)b_j\vc-\sum_{j=1}^3\partial_j(\vn\psi\wedge (b_j\vc))+\sum_{j=1}^3(\partial_j\vn\psi)\wedge (b_j\vc),$$
and Lemma \ref{Lema_VectorialIdentite1} is proven. \hfill $\blacksquare$
%%%%%%%%%%%%%%%%%%%%%%%%%%%%%%%%%%%%%%%%%%%%%%%%%%%%%%%%%%%%%%%%

\end{document}